\documentclass[11pt]{article}
\usepackage[T1]{fontenc}
\usepackage{lmodern}
\usepackage{amsmath,amssymb,amsthm}
\IfFileExists{microtype.sty}{\usepackage{microtype}}{}
\usepackage[margin=1in]{geometry}
\usepackage[hidelinks]{hyperref}
\hypersetup{pdftitle={The Vinogradov system in finite fields and two-source extractors with arbitrarily small min-entropy rate},pdfauthor={Mark Lewko}}

\newtheorem{theorem}{Theorem}
\newtheorem{lemma}[theorem]{Lemma}
\newtheorem{proposition}[theorem]{Proposition}
\newtheorem{corollary}[theorem]{Corollary}
\theoremstyle{remark}

\newcommand{\F}{\mathbb F}

\newcommand{\Z}{\mathbb Z}
\newcommand{\G}{\Gamma}
\newcommand{\eps}{\varepsilon}
\newcommand{\Fbar}{\overline{\F}_p}
\newcommand{\spn}[1]{\langle #1\rangle}
\DeclareMathOperator{\Tr}{Tr}
\DeclareMathOperator{\Ext}{Ext}
\numberwithin{equation}{section}
\allowdisplaybreaks[2]

\title{The Vinogradov system in finite fields and two-source extractors with arbitrarily small min-entropy rate}
\author{Mark Lewko}
\date{}

\begin{document}
\maketitle

\begin{abstract}
Fix $k\ge2$ and $1/(k+1)<\rho\le1$. Let $q=p^n$, where the prime $p$ is fixed and sufficiently large in terms of $k,\rho$, and $n$ tends to infinity through prime values.
We show that for every $A\subseteq\F_q$ with $|A|\ge q^\rho$, the Vinogradov system of degree $k$ in $2s$ variables has at most $|A|^{2s}q^{-k/2-\eta}$ solutions in $A$.
Here the integer $s$ and the saving $\eta>0$ depend only on $k,\rho,p$. At $\rho=1/(k+1)$ an example rules out every saving $\eta>0$ in this moment estimate.
For $k=2$ the author previously proved the analogous estimate over prime fields.
We deduce a power saving for the double sums of $e(xy+(xy)^2+\dots+(xy)^k)$ over any two sets of size at least $q^\rho$, where $e$ is any nontrivial additive character of $\F_q$.
Choosing $k$ sufficiently large gives a polynomial two-source extractor with exponentially small error at every fixed positive min-entropy rate.
The proof of the moment estimate uses the theorem of Gowers, Green, Manners and Tao on Marton's conjecture and a bound for the intersection of the moment curve with an affine subspace over $\F_p$.
\end{abstract}

\section{Introduction}

Let $k\ge1$ and $s\ge1$ be integers.
For a finite subset $A$ of a commutative ring, let $J_{s,k}(A)$ count the solutions of
\begin{equation}\label{eq:system}
 x_1^j+\dots+x_s^j=y_1^j+\dots+y_s^j\qquad(1\le j\le k)
\end{equation}
with $x_1,\dots,x_s,y_1,\dots,y_s\in A$. Thus $k$ is the number of equations and $s$ is the number of summands on each side of an equation.
When $A=\{1,\dots,N\}$, this is the quantity that occurs in Vinogradov's mean value theorem. In this case
\begin{equation}\label{eq:vmvtlower}
 J_{s,k}(A)\gtrsim_{s,k}N^s+N^{2s-k(k+1)/2}.
\end{equation}
The first term follows by taking $y_i=x_i$ for every $i$.
The second term arises from the fact that the sum of the $j$-th powers of $s$ elements of $A$ takes at most $sN^j$ values.
Thus the vector of the first $k$ power sums takes at most $s^kN^{k(k+1)/2}$ values, and the bound follows from the Cauchy--Schwarz inequality.
The main conjecture in Vinogradov's mean value theorem states that \eqref{eq:vmvtlower} is sharp up to a factor of $N^\eps$.
Wooley \cite{WooleyCubic} proved the cubic case, and Bourgain, Demeter and Guth \cite{BDG} proved the case $k\ge4$. See also \cite{WooleyNEC}.

In this paper we study the system \eqref{eq:system} when $A$ is an arbitrary subset of a finite field.
Let $p$ be a prime, let $q=p^n$, and let $\F_q$ denote the field with $q$ elements.
For a nonempty set $A\subseteq\F_q$ the two lower bounds above take the form
\begin{equation}\label{eq:trivial}
 J_{s,k}(A)\ge|A|^s\qquad\text{and}\qquad J_{s,k}(A)\ge|A|^{2s}q^{-k}.
\end{equation}
Here the second bound holds because the vector of power sums takes values in $\F_q^k$.
The quotient $J_{s,k}(A)/|A|^{2s}$ is the probability that \eqref{eq:system} holds when the $2s$ variables are chosen independently and uniformly from $A$.
We will be interested in estimates of the form
\begin{equation}\label{eq:goal}
 J_{s,k}(A)\le|A|^{2s}q^{-k/2-\eta}
\end{equation}
for some $\eta>0$.
For sets of size at least $q^\rho$, we want a single choice of $s$ and $\eta$ that works as $n$ grows.

The exponent $k/2$ arises naturally when we estimate a double exponential sum.
Write $P_k(t):=t+t^2+\dots+t^k$ and $e(z):=\exp(2\pi i\Tr(z)/p)$, where $\Tr$ is the trace from $\F_q$ to $\F_p$ and its value in the exponential is represented by an integer.
Two applications of H\"older's inequality show that for every even integer $s$, all nonempty sets $A,B\subseteq\F_q$ and every nonzero $\lambda\in\F_q$ one has
\begin{equation}\label{eq:bilinear}
 \Bigl|\frac1{|A||B|}\sum_{x\in A}\sum_{y\in B}e\bigl(\lambda P_k(xy)\bigr)\Bigr|^{s^2}\le q^{k/2}\Bigl(\frac{J_{s,k}(A)}{|A|^{2s}}\cdot\frac{J_{s,k}(B)}{|B|^{2s}}\Bigr)^{1/2}
\end{equation}
(see Proposition~\ref{prop:bilinear}).
Thus if \eqref{eq:goal} holds for both $A$ and $B$, then the normalized sum on the left is at most $q^{-\eta/s^2}$.
When $k=1$ we have $J_{2,1}(A)\le|A|^3$, and thus \eqref{eq:goal} holds with $s=2$ for every set with $|A|\ge q^{1/2+\eta}$.
In this case \eqref{eq:bilinear} recovers the classical fact that the sums of $e(\lambda xy)$ over $x\in A$ and $y\in B$ exhibit cancellation when $|A||B|$ exceeds $q$ by a power of $q$.
The threshold $|A|=q^{1/2}$ for $k=1$ already has a simple explanation. If $A$ is a subspace of $\F_q$ over $\F_p$ then $J_{s,1}(A)=|A|^{2s-1}$ for every $s$, and \eqref{eq:goal} fails when $|A|\le q^{1/2}$.

For $k\ge2$ we can use the way the moment curve $\G(x):=(x,x^2,\dots,x^k)$ meets subspaces of $\F_q^k$.
Thus one may hope that \eqref{eq:goal} holds for sets that are much smaller than $q^{1/2}$, at the cost of taking $s$ large.
This idea is due to Bourgain \cite{Bourgain}.
He used the sum-product phenomenon in prime fields to obtain a power saving for the sums in \eqref{eq:bilinear} when $k=2$, $q=p$ is prime, and $|A|,|B|\ge p^{1/2-\delta}$, where $\delta>0$ is an absolute constant.
Estimates for the system \eqref{eq:system} with $k=2$ in arbitrary subsets of $\F_p$ were later obtained by Mansfield and Mudgal \cite{MM}.
In \cite[Lemma~6]{Lewko} the author showed that every nonempty set $A\subseteq\F_p$ satisfies
\begin{equation}\label{eq:quadratic}
 J_{s,2}(A)\lesssim_s|A|^{2s-3+2^{2-s}}+|A|^{2s-1}p^{-1}
\end{equation}
for every $s\ge2$.
This is a consequence of a sharp form of the Szemer\'edi--Trotter theorem in $\F_p^2$.
In particular, when $k=2$ and $q=p$ the estimate \eqref{eq:goal} holds, with $s$ and $\eta$ depending on $\rho$, for all sets with $|A|\ge p^\rho$ whenever $\rho>\frac13$.

An interval in a prime field explains why the moment estimate requires $\rho>1/(k+1)$, which for $k=2$ means $\rho>1/3$.
If $A=\{1,\dots,N\}\subseteq\F_p$, every integer solution of \eqref{eq:system} gives a solution modulo $p$.
Thus \eqref{eq:vmvtlower} gives $J_{s,k}(A)/|A|^{2s}\gtrsim_{s,k}N^{-D}$, where $D:=k(k+1)/2$.
At $N=\lceil p^{1/(k+1)}\rceil$ this is $\gtrsim_{s,k}p^{-k/2}$, for every fixed $s$.
We obtain the corresponding upper bound for every $\rho>1/(k+1)$ in fields of large prime degree over $\F_p$, with $p$ fixed and sufficiently large in terms of $k$ and $\rho$.

\begin{theorem}\label{thm:main}
Let $k\ge2$ be an integer, let $\frac1{k+1}<\rho\le1$, and let $p$ be a prime such that
\begin{equation}\label{eq:plarge}
 p^{\rho-\frac1{k+1}}>24k+12.
\end{equation}
Then there exist an even integer $s$ and a real number $\eta>0$, depending only on $k$, $\rho$ and $p$, such that the following holds.
If $n$ is a sufficiently large prime, $q=p^n$, and $A\subseteq\F_q$ satisfies $|A|\ge q^\rho$, then
\[
 J_{s,k}(A)\le|A|^{2s}q^{-k/2-\eta}.
\]
\end{theorem}

For instance, if $\rho=1/k$ then the theorem applies to every prime $p>(24k+12)^{k(k+1)}$.
The proof shows that one may take
\[
 \eta=\frac{k^2+k-1}4\min\Bigl\{\gamma,\frac1{k(k+1)}\Bigr\},
\]
where $\gamma>0$ is defined by $p^\gamma=p^{\rho-\frac1{k+1}}/(24k+12)$.
In particular, if $\rho=1/k$ and $p>(24k+12)^{4k(k+1)}$ then one may take $\eta=\frac18$.

The role of intervals in these extensions is played by polynomials of small degree in a generator of the field.
Let $\theta$ generate $\F_q$ over $\F_p$. For $1\le h\le n$, let $I_h$ denote the set of elements $a_0+a_1\theta+\dots+a_{h-1}\theta^{h-1}$ with $a_i\in\F_p$.
The sum of the $j$-th powers of $s$ elements of $I_h$ is a polynomial in $\theta$ of degree at most $j(h-1)$ with coefficients in $\F_p$.
There are at most $p^{j(h-1)+1}$ possibilities for the $j$-th sum and at most $p^{D(h-1)+k}$ for the whole vector. Cauchy--Schwarz gives
\begin{equation}\label{eq:example}
 J_{s,k}(I_h)\ge|I_h|^{2s}p^{-D(h-1)-k}
\end{equation}
for every $s$.
Take $h=\lceil n/(k+1)\rceil$. Then $|I_h|\ge q^{1/(k+1)}$ and $D(h-1)+k<kn/2+k$.
Thus $J_{s,k}(I_h)/|I_h|^{2s}\ge p^{-k}q^{-k/2}$ for every $s$.
With $p$ fixed, this rules out any saving $q^{-\eta}$ at the endpoint itself.

Proper subfields give a different obstruction.
Indeed, if $d$ divides $n$ and $A$ is the subfield of $\F_q$ with $q^{1/d}$ elements, then the power sums of elements of $A$ lie in $A$, which gives $J_{s,k}(A)\ge|A|^{2s}q^{-k/d}$ for every $s$.
Thus the conclusion of the theorem fails whenever $n$ has a divisor $d$ with $2\le d\le1/\rho$.

For $p>k$ we do not know whether the restriction \eqref{eq:plarge} on the characteristic is necessary.
It arises from a factor of $(24k+12)^n$ in Theorem~\ref{thm:section} below, which in turn comes from a bound for the number of zeros and poles of the zeta function of an algebraic variety.
The value of $s$ given by the proof is very large.
The proof gives $s\le\exp(O(k^3p^3/\eta))$. A natural question is whether $D=k(k+1)/2$ summands on each side suffice for the moment estimate for every $\rho>1/(k+1)$ (see Section~\ref{sec:remarks}).

The exponential sum estimate now follows from \eqref{eq:bilinear}.

\begin{theorem}\label{thm:sums}
Let $k$, $\rho$, $p$, $s$ and $\eta$ be as in Theorem~\ref{thm:main}, let $n$ be a sufficiently large prime and let $q=p^n$.
If $A,B\subseteq\F_q$ satisfy $|A|\ge q^\rho$ and $|B|\ge q^\rho$, then for every nonzero $\lambda\in\F_q$ we have
\[
 \Bigl|\sum_{x\in A}\sum_{y\in B}e\bigl(\lambda P_k(xy)\bigr)\Bigr|\le|A||B|\,q^{-\eta/s^2}.
\]
\end{theorem}

In particular, for each fixed $\rho>0$ we may choose $k$ with $\frac1{k+1}<\rho$ and obtain a power saving for all sets of size at least $q^\rho$, provided that $p$ is sufficiently large in terms of $\rho$.
Estimates of this type are closely related to the problem of constructing two-source extractors in theoretical computer science.
For a random variable $X$ on a finite set $\Omega$, its min-entropy is $H_\infty(X):=-\log_2\max_x\Pr(X=x)$.
A function $\Ext:\Omega\times\Omega\to\Omega'$ is a two-source extractor for min-entropy $h$ with error $\eps$ if the distribution of $\Ext(X,Y)$ is within $\eps$ of uniform in total variation distance whenever $X,Y$ are independent and each has min-entropy at least $h$.
The min-entropy rate is $h/\log_2|\Omega|$. For a source uniform on $A\subseteq\F_q$, rate at least $\rho$ means simply that $|A|\ge q^\rho$.

Chor and Goldreich \cite{CG} showed that the inner product extracts a bit with exponentially small error from two sources of any fixed rate greater than $1/2$. Bourgain \cite{Bourgain} was the first to go below $1/2$, using exponential sums over prime fields. Chattopadhyay and Zuckerman \cite{CZ} later gave an explicit extractor for polylogarithmic min-entropy with error polynomially small in the input length.
Returning to Bourgain's approach, the author obtained every rate above $21/44$ for the paraboloid in $\F_p^3$ when $p\equiv3\pmod4$ \cite{Lewko49}. The $4/9$ claim in that paper for the paraboloid in $\F_p^4$ used an incorrect energy estimate of Rudnev and Shkredov \cite[Lemma~4]{RS}. Before the author's Szemer\'edi--Trotter theorem \cite{Lewko}, the best published min-entropy threshold known to the author for an explicit two-source extractor with exponentially small error was $123/260$, obtained from the rectangle estimate in \cite{LewkoRect} for the three-dimensional paraboloid. The Szemer\'edi--Trotter theorem gave every rate above $3/8$ for that paraboloid and every rate above $1/3$ for Bourgain's scalar quadratic $xy+x^2y^2$ over every odd prime field \cite[Remark~7 and Corollary~5]{Lewko}.
Chattopadhyay and Gurumukhani \cite{ChG} obtained exponentially small error at every rate above $3/7$ for the cubic $P_3(xy)$ and conjectured that higher degrees suffice at every fixed positive rate over prime fields. Corollary~\ref{cor:extractor} gives the corresponding conclusion over the extensions considered here, using a linear projection for the output.

\begin{corollary}\label{cor:extractor}
Let $k,\rho,p,s,\eta,n,q$ be as in Theorem~\ref{thm:sums}, and let $\pi:\F_q\to\F_p^m$ be a surjective $\F_p$-linear map, where $1\le m\le n$. Then
\[
 \Ext(x,y):=\pi\bigl(xy+(xy)^2+\dots+(xy)^k\bigr)
\]
is a two-source extractor for min-entropy rate $\rho$ with error at most $\frac12p^{m/2}q^{-\eta/s^2}$.
\end{corollary}

For $n\ge s^2/\eta$ we may take $m=\lfloor\eta n/s^2\rfloor\ge1$. The error is then at most $\frac12q^{-\eta/(2s^2)}$, which is exponentially small in the input length $n\log_2p$.
The ratio of the output length to the input length is $m/n$, which tends to $\eta/s^2$ as $n$ grows.

The example $I_h$ settles the threshold for the moment estimate. The threshold for extraction by $P_2(xy)$ is still open. Section~\ref{sec:lower} shows that the worst error at rate $\rho$ is at least a constant times $q^{1/4-\rho}$ for $1/4\le\rho<1$. This holds for every nonzero linear coordinate when $p$ is fixed, and for the half-interval and parity bits when $q=p$ tends to infinity. The constant may depend on $p$ in the first case and is absolute in the second. In particular each construction has constant error at rate $1/4$. We do not know whether either gives a power saving at every fixed rate above $1/4$.

While this paper was being finalized, Li \cite{Li} gave a different explicit construction for binary sources at every fixed positive entropy rate, with linear output length and exponentially small error. His output length is asymptotically optimal in terms of the source entropy and the error. His proof also uses the characteristic-two case of Marton's conjecture.
Li also reports independent work of Doron and Goodman for entropy $h\ge N^{1-\alpha}$ on $N$-bit inputs, with error $2^{-h^\alpha}$ for some $\alpha>0$.
An earlier connection between Marton's conjecture and two-source extraction was made by Ben-Sasson and Ron-Zewi \cite{BR} through approximate duality. See also \cite{AGMR}.

To prove Theorem~\ref{thm:main}, let $\mu$ be the distribution of $(x,x^2,\dots,x^k)$ for $x$ uniform on $A$. Write $\mu^{*s}$ for the distribution of the sum of $s$ independent points with this law.
Then $J_{s,k}(A)/|A|^{2s}$ is the sum of the squares of the point masses of $\mu^{*s}$, which we denote by $\|\mu^{*s}\|_2^2$.
This quantity does not increase with $s$.
For a fixed set $A$ it converges, as $s\to\infty$, to $1/|W|$, where $W$ is the subspace of $\F_q^k$ over $\F_p$ spanned by the differences of the points of the support of $\mu$.
Indeed, the characters whose average against $\mu$ has absolute value one are precisely those constant on its support. The other Fourier coefficients tend to zero under convolution. We need the collision probability to fall below $q^{-k/2}$ after a number of steps independent of $n$.
Suppose doubling the number of summands reduces the collision probability by at most a factor $L$, so that $\|\mu^{*2j}\|_2^2\ge L^{-1}\|\mu^{*j}\|_2^2$.
By a lemma of Bourgain and Gamburd, which is based on the Balog--Szemer\'edi--Gowers theorem, a large part of the measure $\mu^{*j}$ is then carried by a set of small doubling.
Bourgain used sum-product estimates at this point in the prime-field argument.
Here $\F_q^k$ is a vector space over the fixed field $\F_p$, so we can use the theorem of Gowers, Green, Manners and Tao on Marton's conjecture \cite{GGMT2}. It says that a set with doubling constant $L$ is covered by $L^{O_p(1)}$ cosets of a subspace of no larger size.
The polynomial dependence on $L$ lets us take $L$ to be a small power of $q$ and bound the number of doublings in terms of $k,\rho,p$ alone.
It follows that some affine subspace $V$ has $|V|\le L^{O_p(1)}\|\mu^{*j}\|_2^{-2}$ and $\mu(V)\ge L^{-O_p(1)}$.
This reduces Theorem~\ref{thm:main} to the following estimate, which is the main new ingredient in this paper.

\begin{theorem}\label{thm:section}
Let $k\ge2$ be an integer, let $p>k$ and $n$ be primes, and let $q=p^n$.
Let $V\subseteq\F_q^k$ be an affine subspace over $\F_p$ of dimension $r<(k+1)n/2$.
Then the number of $x\in\F_q$ such that $(x,x^2,\dots,x^k)\in V$ is at most
\[
 3(24k+12)^n\,p^{r/D+1}.
\]
\end{theorem}

Sets that have few points in every affine subspace of a given dimension are called subspace-evasive (see \cite{PR}, \cite{AGMR}).
Over $\F_2$ it is known that a suitable set of this type gives rise to a two-source extractor with exponentially small error, by the results of \cite{BR} and \cite{GGMT1} (see \cite[Theorem~10]{AGMR}).
The argument outlined above is a variant of this principle that applies to the system \eqref{eq:system} and in odd characteristic.

Apart from the factor $3(24k+12)^np$ the estimate in Theorem~\ref{thm:section} is sharp.
Indeed, let $I_h$ be the set defined above and suppose that $k(h-1)\le n-2$.
Then the $j$-th powers of the elements of $I_h$ lie in $I_{j(h-1)+1}$.
Thus the image of $I_h$ lies in a subspace of dimension $r=D(h-1)+k<(k+1)n/2$, while $|I_h|=p^h>p^{r/D}$.

To see where the exponent comes from, suppose $U\subseteq\F_q$ is a subspace over $\F_p$ of dimension $h$ such that $(u,u^2,\dots,u^k)\in V$ for all $u\in U$.
Then $V$ is a linear subspace.
If $u\in U$ then $tu\in U$ for every $t\in\F_p$, and thus $V$ contains the vectors $(tu,t^2u^2,\dots,t^ku^k)$.
Since $p>k$ we may interpolate in $t$, and we find that $V$ contains the vector whose $j$-th coordinate is $u^j$ and whose other coordinates vanish.
Thus $V$ contains the product of the spaces $U_j$, $1\le j\le k$, where $U_j$ is the span of the $j$-th powers of the elements of $U$.
Since $p>k$ this is also the span of the products of $j$ elements of $U$.
A linear analogue of the Cauchy--Davenport inequality, which is due to Hou, Leung and Xiang \cite{HLX} and which uses that $n$ is prime, shows that $\dim U_j\ge\min\{n,1+j(h-1)\}$.
Summing over $j$ we obtain $r\ge D(h-1)+k$ whenever $1+k(h-1)\le n$, that is, $|U|\le p^{(r-k)/D+1}$.

For a general set we use differentiation in place of dilation. Writing $x$ in a basis of $\F_q$ over $\F_p$ turns the condition $\G(x)\in V$ into a system of polynomial equations in $n$ variables.
At a generic point of an irreducible set of solutions of dimension $h$, the first derivatives of $x$ span a space of dimension $h$. A triangular change of coordinates expresses translation on the moment curve. In these coordinates the $j$-th derivative has $j!u^j$ in coordinate $j$ and zero in every later coordinate, where $u$ runs through that $h$-dimensional space.
We can then count dimensions one coordinate at a time, just as above.
There are two points to address in making this argument precise.
The product estimate requires a field extension of degree $n$. After passing to an algebraic closure, the polynomial defining $\F_q$ splits.
We work instead over the union of the fields $\F_{p^m}$ with $n\nmid m$, over which this polynomial stays irreducible.
We then need a point count uniform in $n$. Rationality of the zeta function allows us to pass from estimates over these auxiliary fields to an estimate over $\F_p$ with a constant given by a theorem of Wan and Zhang \cite{WZ}. This gives the factor $3(24k+12)^n$.

Sections~\ref{sec:flat} and~\ref{sec:proofs} deduce the moment and extractor estimates from Theorem~\ref{thm:section}.
We prove that theorem in Sections~\ref{sec:dimension} and~\ref{sec:counting}. Sections~\ref{sec:lowerfixed} and~\ref{sec:lowerprime} give lower bounds for the extractor error in fixed characteristic and over prime fields, respectively. Section~\ref{sec:remarks} discusses the choice of $s$ and the restriction on $p$.
Appendix~\ref{app:lowerproofs} proves the auxiliary results used in the scalar lower bounds and determines when a linear coordinate can be constant on two affine sources. Appendix~\ref{app:related} treats the cubic pairing on the parabola and the examples for the paraboloid.

Subspaces of $\F_q^k$ are taken over $\F_p$ unless another field is specified.
For subspaces $U,U'$ of an extension of a field $F$, let $\spn{UU'}$ denote the $F$-linear span of the products $uu'$ with $u\in U$ and $u'\in U'$. We use $\spn{U^j}$ for the span of products of $j$ elements of $U$.
In the notation $X\lesssim Y$, the implicit constant is absolute unless dependence on parameters is indicated by subscripts.

\subsection*{Acknowledgments and AI usage}
This paper grew out of discussions with OpenAI's ChatGPT models. The author used Anthropic's Claude to check the arguments and improve the exposition. The author takes responsibility for all content.

\section{A flattening lemma}\label{sec:flat}

Let $G$ be a finite abelian group.
By a probability measure on $G$ we mean a function $\mu:G\to[0,1]$ with $\sum_{x\in G}\mu(x)=1$, and we write $\mu(E):=\sum_{x\in E}\mu(x)$ for $E\subseteq G$.
For functions $\mu$ and $\nu$ on $G$ we define
\[
 (\mu*\nu)(x):=\sum_{y\in G}\mu(y)\nu(x-y),\qquad\tilde\mu(x):=\mu(-x),\qquad\|\mu\|_2:=\Bigl(\sum_{x\in G}|\mu(x)|^2\Bigr)^{1/2},
\]
and we let $\mu^{*s}$ denote the convolution of $s$ copies of $\mu$.
If $X$ and $Y$ are independent random variables with distributions $\mu$ and $\nu$, then $\mu*\nu$ is the distribution of $X+Y$, and $\|\mu\|_2^2$ is the probability that two independent copies of $X$ coincide.
We will use the following two elementary facts.
If $\mu$ and $\nu$ are probability measures on $G$ and $E\subseteq G$, then
\begin{equation}\label{eq:young}
 \|\mu*\nu\|_2\le\|\nu\|_2\qquad\text{and}\qquad(\mu*\nu)(E)\le\max_{z\in G}\nu(z+E).
\end{equation}
Indeed, $\mu*\nu$ is an average of translates of $\nu$.
In particular, $\|\mu^{*s}\|_2$ does not increase with $s$.

The following lemma of Bourgain and Gamburd describes what happens when doubling the number of summands gives only a small decrease.
It first appeared implicitly in \cite{BG}, and its proof is based on the Balog--Szemer\'edi--Gowers theorem.
We state it for abelian groups, in the form given in \cite[Lemma~3.8]{Breuillard} (see \cite[Lemma~15]{Varju} for a proof).

\begin{theorem}[Bourgain and Gamburd]\label{thm:BG}
There exists an absolute constant $C_0\ge1$ such that the following holds.
Let $\nu$ be a probability measure on a finite abelian group $G$ and let $T\ge2$ be such that $\|\nu*\nu\|_2\ge T^{-1}\|\nu\|_2$.
Then there exists a nonempty set $S\subseteq G$ such that
\[
 |S|\le T^{C_0}\|\nu\|_2^{-2},\qquad|S+S|\le T^{C_0}|S|\qquad\text{and}\qquad(\nu*\tilde\nu)(x)\ge\frac1{T^{C_0}|S|}\quad\text{for all }x\in S.
\]
\end{theorem}

We use Marton's conjecture in the bounded-torsion form proved by Gowers, Green, Manners and Tao \cite[Theorem~1.1]{GGMT2}. Their earlier paper \cite{GGMT1} treats groups of exponent two.

\begin{theorem}[Gowers, Green, Manners and Tao]\label{thm:Marton}
Let $d\ge2$ be an integer and let $G$ be an abelian group such that $dx=0$ for every $x\in G$.
Suppose that $S\subseteq G$ is a finite nonempty set with $|S+S|\le T|S|$.
Then $S$ is covered by at most $(2T)^{cd^3}$ cosets of a subgroup $H$ of $G$ with $|H|\le|S|$, where $c$ is an absolute constant.
\end{theorem}

We record the consequence for a probability measure on a vector space over a fixed prime field.

\begin{proposition}\label{prop:inverse}
For every prime $p$ there exists a constant $C_p\ge1$ such that the following holds.
Let $G$ be a finite-dimensional vector space over $\F_p$, let $\nu$ be a probability measure on $G$, and let $L\ge2$ be such that
\[
 \|\nu*\nu\|_2^2\ge L^{-1}\|\nu\|_2^2.
\]
Then there exist a subspace $H$ of $G$ and an element $x\in G$ such that
\[
 |H|\le L^{C_p}\|\nu\|_2^{-2}\qquad\text{and}\qquad\nu(x+H)\ge L^{-C_p}.
\]
Moreover, one may take $C_p\lesssim p^3$.
\end{proposition}

\begin{proof}
The hypothesis gives $\|\nu*\nu\|_2\ge L^{-1/2}\|\nu\|_2$, so Theorem~\ref{thm:BG} applies with $T=L$.
Let $S$ be the set given by that theorem.
Summing the third property over $x\in S$ we see that $(\nu*\tilde\nu)(S)\ge L^{-C_0}$.
By Theorem~\ref{thm:Marton}, applied with $d=p$ and $T=L^{C_0}$, the set $S$ is covered by at most $(2L^{C_0})^{cp^3}\le L^{(C_0+1)cp^3}$ cosets of a subgroup $H$ of $G$ with $|H|\le|S|\le L^{C_0}\|\nu\|_2^{-2}$.
Every additive subgroup of $G$ is closed under multiplication by scalars in $\F_p$, so $H$ is a subspace.
For one of these cosets, say $y+H$, we have
\[
 (\nu*\tilde\nu)(y+H)\ge L^{-C_0-(C_0+1)cp^3}.
\]
By the second part of \eqref{eq:young}, applied to $\tilde\nu*\nu=\nu*\tilde\nu$, there is an $x\in G$ with $\nu(x+H)\ge(\nu*\tilde\nu)(y+H)$.
This gives the stated bound with $C_p:=C_0+(C_0+1)cp^3$.
\end{proof}

We now apply the proposition at successive powers of two. A failure of the collision probability to decrease forces concentration on a translate of a subspace.

\begin{lemma}\label{lem:flat}
Let $p$ be a prime and let $C_p$ be the constant in Proposition~\ref{prop:inverse}.
Let $G$ be a finite-dimensional vector space over $\F_p$, let $\mu$ be a probability measure on $G$, and let $M\ge1$ and $L\ge2$ be real numbers.
Suppose that
\begin{equation}\label{eq:nonconc}
 \mu(x+H)<L^{-C_p}
\end{equation}
for every $x\in G$ and every subspace $H$ of $G$ with $|H|\le M$.
Then for every integer $t\ge1$ we have
\[
 \|\mu^{*2^t}\|_2^2\le\max\bigl\{L^{-t},\,L^{C_p}M^{-1}\bigr\}.
\]
\end{lemma}

\begin{proof}
Let $R_j:=\|\mu^{*j}\|_2^{-2}$.
By \eqref{eq:young} the sequence $R_j$ is non-decreasing, and $R_1\ge1$.
Suppose first that $R_{2j}>LR_j$ for every $j=2^i$ with $0\le i<t$.
Then $R_{2^t}>L^tR_1\ge L^t$.
Otherwise there is a $j=2^i$ with $0\le i<t$ such that $R_{2j}\le LR_j$.
We apply Proposition~\ref{prop:inverse} to the measure $\nu:=\mu^{*j}$, and obtain a subspace $H$ and an element $x$ with $|H|\le L^{C_p}R_j$ and $\nu(x+H)\ge L^{-C_p}$.
By the second part of \eqref{eq:young} there is a $z\in G$ with $\mu(z+H)\ge\nu(x+H)\ge L^{-C_p}$.
In view of \eqref{eq:nonconc} this implies that $|H|>M$.
Thus $R_{2^t}\ge R_j\ge L^{-C_p}|H|>L^{-C_p}M$.
\end{proof}

\section{Proof of the main theorems}\label{sec:proofs}

We prove the moment estimate directly for the probability measures that occur in the extractor statement.
For a probability measure $\mu_0$ on $\F_q$, write $\G_*\mu_0$ for the measure on $\F_q^k$ that puts mass $\mu_0(x)$ at $\G(x)$ for every $x$.
If $\mu_0$ is the uniform probability measure on a nonempty set $A\subseteq\F_q$ and $\mu=\G_*\mu_0$, then
\begin{equation}\label{eq:Jnorm}
 J_{s,k}(A)=|A|^{2s}\|\mu^{*s}\|_2^2.
\end{equation}

Let $k$, $\rho$ and $p$ be as in Theorem~\ref{thm:main}, and let $C_p$ be the constant in Proposition~\ref{prop:inverse}.
We define $\gamma>0$ by $p^\gamma=p^{\rho-\frac1{k+1}}/(24k+12)$, and we let
\begin{equation}\label{eq:parameters}
 \tau:=\frac12\min\Bigl\{\gamma,\frac1{k(k+1)}\Bigr\},\quad
 \eta:=\Bigl(D-\frac12\Bigr)\tau,\quad
 t:=\Bigl\lceil\frac{2(k+1)C_p}\tau\Bigr\rceil,\quad s:=2^t.
\end{equation}
Note that $p>24k+12>k$ by \eqref{eq:plarge}.

\begin{proposition}\label{prop:measure}
Let $n$ be a prime, let $q=p^n$, and suppose that $q^{\tau/2}>3p$ and $q^{\tau/(2C_p)}\ge2$.
Let $\mu_0$ be a probability measure on $\F_q$ such that $\mu_0(x)\le q^{-\rho}$ for every $x\in\F_q$, and let $\mu:=\G_*\mu_0$.
Then
\[
 \|\mu^{*s}\|_2^2\le q^{-k/2-\eta}.
\]
\end{proposition}

\begin{proof}
Let $a:=\frac k2+D\tau$.
Since $D\tau\le\frac14$ we have $a<\frac{k+1}2$.
Let $H$ be a subspace of $\F_q^k$ with $|H|\le q^a$ and let $v\in\F_q^k$.
The dimension $r$ of $H$ satisfies $r\le an<(k+1)n/2$, and thus Theorem~\ref{thm:section} applies to the affine subspace $v+H$.
We have $a/D=1/(k+1)+\tau$ and $(24k+12)^n=q^{\rho-1/(k+1)-\gamma}$. Since each point has $\mu_0$-mass at most $q^{-\rho}$, Theorem~\ref{thm:section} gives
\[
 \mu(v+H)\le 3p\,q^{-\rho}(24k+12)^n q^{a/D}
          =3p\,q^{-\gamma+\tau}\le3p\,q^{-\tau},
\]
where the last inequality uses $\tau\le\gamma/2$.
Let $L:=q^{\tau/(2C_p)}$.
By our assumptions on $q$ we have $L\ge2$ and $\mu(v+H)\le3p\,q^{-\tau}<q^{-\tau/2}=L^{-C_p}$.
Thus the hypothesis of Lemma~\ref{lem:flat} is satisfied with $G=\F_q^k$ and $M=q^a$, and we obtain
\[
 \|\mu^{*s}\|_2^2\le\max\bigl\{q^{-t\tau/(2C_p)},\,q^{\tau/2-a}\bigr\}.
\]
By the choice of $t$ we have $t\tau/(2C_p)\ge k+1>a$.
Moreover, $a-\tau/2=\frac k2+(D-\frac12)\tau=\frac k2+\eta$.
Thus $\|\mu^{*s}\|_2^2\le q^{-k/2-\eta}$.
\end{proof}

\begin{proof}[Proof of Theorem~\ref{thm:main}]
Let $n$ be a prime that is large enough for the assumptions of Proposition~\ref{prop:measure} to hold, and let $A\subseteq\F_q$ with $|A|\ge q^\rho$.
The uniform probability measure $\mu_0$ on $A$ satisfies $\mu_0(x)\le q^{-\rho}$ for every $x$.
The theorem now follows from Proposition~\ref{prop:measure} and \eqref{eq:Jnorm}.
\end{proof}

We now turn to exponential sums.
For $u,v\in\F_q^k$ we write $u\cdot v:=u_1v_1+\dots+u_kv_k$, so that
\[
 \G(x)\cdot\G(y)=\sum_{j=1}^k(xy)^j.
\]
Taking the two measures in the next proposition to be uniform on $A$ and $B$ gives \eqref{eq:bilinear}. We include the argument from \cite{Bourgain} (see also \cite[Section~6.3]{Lewko}).

\begin{proposition}\label{prop:bilinear}
Let $\alpha_0$ and $\beta_0$ be probability measures on $\F_q$, and let $\alpha:=\G_*\alpha_0$ and $\beta:=\G_*\beta_0$.
Then for every even integer $s\ge2$ and every nonzero $\lambda\in\F_q$ we have
\[
 \Bigl|\sum_{x,y\in\F_q}\alpha_0(x)\beta_0(y)\,e\bigl(\lambda P_k(xy)\bigr)\Bigr|^{s^2}\le q^{k/2}\,\|\alpha^{*s}\|_2\,\|\beta^{*s}\|_2.
\]
\end{proposition}

\begin{proof}
Let $\Sigma$ denote the sum on the left.
Then $\Sigma=\sum_{u,v\in\F_q^k}\alpha(u)\beta(v)e(\lambda u\cdot v)$.
Let $\alpha'$ denote the distribution of $u_1+\dots+u_{s/2}-u_{s/2+1}-\dots-u_s$, where $u_1,\dots,u_s$ are independent random variables with distribution $\alpha$, and define $\beta'$ similarly.
Expanding the $s$-th power of the absolute value, we see that
\begin{equation}\label{eq:expand}
 \Bigl|\sum_v\beta(v)e(\lambda u\cdot v)\Bigr|^s=\sum_v\beta'(v)e(\lambda u\cdot v)
\end{equation}
for every $u\in\F_q^k$, and the same holds with the roles of $\alpha$ and $\beta$ interchanged.
By H\"older's inequality and \eqref{eq:expand} we have
\[
 |\Sigma|^s\le\sum_u\alpha(u)\Bigl|\sum_v\beta(v)e(\lambda u\cdot v)\Bigr|^s=\sum_v\beta'(v)\sum_u\alpha(u)e(\lambda u\cdot v).
\]
Applying H\"older's inequality and \eqref{eq:expand} once more, we obtain
\[
 |\Sigma|^{s^2}\le\sum_v\beta'(v)\Bigl|\sum_u\alpha(u)e(\lambda u\cdot v)\Bigr|^s=\sum_{u,v}\alpha'(u)\beta'(v)e(\lambda u\cdot v).
\]
Since $\lambda\ne0$, the orthogonality of the characters of $\F_q^k$ gives
\[
 \sum_u\Bigl|\sum_v\beta'(v)e(\lambda u\cdot v)\Bigr|^2=q^k\|\beta'\|_2^2.
\]
Thus the Cauchy--Schwarz inequality yields $|\Sigma|^{s^2}\le q^{k/2}\|\alpha'\|_2\|\beta'\|_2$.
Finally, $\|\alpha'\|_2^2$ is the probability that two independent random variables with distribution $\alpha'$ coincide.
Moving the terms with a negative sign to the other side of this equation, we see that $\|\alpha'\|_2=\|\alpha^{*s}\|_2$, and similarly $\|\beta'\|_2=\|\beta^{*s}\|_2$.
This completes the proof.
\end{proof}

\begin{proof}[Proof of Theorem~\ref{thm:sums} and Corollary~\ref{cor:extractor}]
Let $n$ be as in the proof of Theorem~\ref{thm:main}, and let $X$ and $Y$ be independent random variables in $\F_q$ whose distributions $\alpha_0$ and $\beta_0$ satisfy $\alpha_0(x)\le q^{-\rho}$ and $\beta_0(x)\le q^{-\rho}$ for every $x$.
By Proposition~\ref{prop:bilinear} and Proposition~\ref{prop:measure}, for every nonzero $\lambda\in\F_q$ we have
\begin{equation}\label{eq:bias}
 \Bigl|\sum_{x,y\in\F_q}\alpha_0(x)\beta_0(y)\,e\bigl(\lambda P_k(xy)\bigr)\Bigr|\le\bigl(q^{k/2}q^{-k/2-\eta}\bigr)^{1/s^2}=q^{-\eta/s^2}.
\end{equation}
Theorem~\ref{thm:sums} is the special case in which $\alpha_0$ and $\beta_0$ are the uniform probability measures on $A$ and $B$.

Now let $\pi:\F_q\to\F_p^m$ be as in Corollary~\ref{cor:extractor}, and let $\omega$ denote the distribution of $\pi(P_k(XY))$.
Let $\chi$ be a nontrivial character of $\F_p^m$.
Then $\chi\circ\pi$ is a nontrivial additive character of $\F_q$, since $\pi$ is surjective and linear over $\F_p$.
Thus $\chi\circ\pi$ is of the form $z\mapsto e(\lambda z)$ for some nonzero $\lambda\in\F_q$, and \eqref{eq:bias} shows that $|\sum_w\omega(w)\chi(w)|\le q^{-\eta/s^2}$.
By the Cauchy--Schwarz inequality and Parseval's identity we have
\[
 \sum_w|\omega(w)-p^{-m}|
 \le p^{m/2}\|\omega-p^{-m}\|_2
 =\Bigl(\sum_{\chi\ne1}\Bigl|\sum_w\omega(w)\chi(w)\Bigr|^2\Bigr)^{1/2}
 \le p^{m/2}q^{-\eta/s^2}.
\]
The total variation distance is half of the left-hand side.
\end{proof}

\section{A dimension estimate}\label{sec:dimension}

In this section and the next we prove Theorem~\ref{thm:section}.
We fix an integer $k\ge2$, primes $p>k$ and $n$, and we let $q=p^n$.
We also fix an element $\theta\in\F_q$ that does not belong to $\F_p$.
Since $n$ is prime, $\theta$ generates $\F_q$ over $\F_p$.
Let $f\in\F_p[T]$ denote its minimal polynomial, which is irreducible of degree $n$.
Finally, we fix an affine subspace $V=v+V_0$ of $\F_q^k$, where $V_0$ is a subspace of dimension $r$.

Given a field $F$ that contains $\F_p$, we let $K_F:=F[T]/(f)$ and we continue to write $\theta$ for the image of $T$ in $K_F$.
Thus $K_F$ is a ring that contains $F$ and $\F_q=K_{\F_p}$, and the elements $1,\theta,\dots,\theta^{n-1}$ form a basis of $K_F$ over $F$.
It is a field if and only if $f$ is irreducible over $F$.
A basis of $\F_q^k$ over $\F_p$ is also a basis of $K_F^k$ over $F$.
We let $FV_0$ denote the span of $V_0$ over $F$ in $K_F^k$, which has dimension $r$ over $F$, and we let $V_F:=v+FV_0$.
Every linear functional $\ell:\F_q^k\to\F_p$ extends uniquely to an $F$-linear functional $\ell_F:K_F^k\to F$, and a vector $w\in K_F^k$ belongs to $V_F$ if and only if
\begin{equation}\label{eq:functionals}
 \ell_F(w)=\ell(v)\qquad\text{for every linear functional }\ell:\F_q^k\to\F_p\text{ that vanishes on }V_0.
\end{equation}
For $z=(z_0,\dots,z_{n-1})\in F^n$ we write
\[
 X(z):=z_0+z_1\theta+\dots+z_{n-1}\theta^{n-1}\in K_F.
\]
The coordinates of the powers $X(z)^j$ in the basis $1,\theta,\dots,\theta^{n-1}$ are polynomials of degree $j$ in $z_0,\dots,z_{n-1}$ with coefficients in $\F_p$ that do not depend on $F$.
Thus, by \eqref{eq:functionals}, the condition
\begin{equation}\label{eq:condition}
 \G(X(z))\in V_F
\end{equation}
is equivalent to a system of polynomial equations in $z_0,\dots,z_{n-1}$ of degree at most $k$ with coefficients in $\F_p$.
We may take this system to consist of $kn-r$ equations.
The quantity to be estimated in Theorem~\ref{thm:section} is the number of solutions of this system in $\F_p^n$.

Proposition~\ref{prop:dimension} bounds the number of algebraically independent coordinates of a solution of \eqref{eq:condition}. When $r<(k+1)n/2$ this number is at most $r/D+1$, provided that $f$ remains irreducible over the field generated by the coordinates over the perfect base field.
We will need two lemmas.
The first is a special case of the theorem of Hou, Leung and Xiang \cite[Theorem~2.4]{HLX}, which is a linear analogue of Kneser's theorem (see also \cite{EL}, \cite{BSZ}).
In the case of extensions of prime degree it is an analogue of the Cauchy--Davenport inequality, and it has the following short proof.

\begin{lemma}\label{lem:products}
Let $E$ be a field extension of a field $F$ of prime degree $n$, and let $U$ and $U'$ be nonzero subspaces of $E$ over $F$.
Then
\[
 \dim_F\spn{UU'}\ge\min\{n,\dim_FU+\dim_FU'-1\}.
\]
In particular, if $\dim_FU=h\ge1$ then $\dim_F\spn{U^j}\ge\min\{n,1+j(h-1)\}$ for every $j\ge1$.
\end{lemma}

\begin{proof}
We argue by induction on $\dim_FU'$.
If $\dim_FU'=1$ then $\spn{UU'}=u'U$ for a nonzero $u'\in U'$, and the claim holds.
If $U=E$ then $\spn{UU'}=E$, and the claim holds as well.
Suppose that $\dim_FU'\ge2$ and $U\ne E$.
We claim that there exist nonzero elements $u\in U$ and $u'\in U'$ such that $e:=u/u'$ satisfies $eU'\not\subseteq U$.
Indeed, suppose that this fails and fix a nonzero $u'\in U'$.
Then $(w/u')U\subseteq U$ for every $w\in U'$.
The set of $t\in E$ with $tU\subseteq U$ is a subring containing $F$. It is a field, since multiplication by a nonzero such $t$ maps the finite-dimensional space $U$ injectively onto itself.
It contains the space $\{w/u':w\in U'\}$, which has dimension at least two.
Since $n$ is prime this field must be equal to $E$.
Thus $EU\subseteq U$, which contradicts the fact that $0\ne U\ne E$.
This proves the claim.

Now let $U_1:=U+eU'$ and $U_1':=U'\cap e^{-1}U$.
Then $U_1'$ contains $u'$, and it is a proper subspace of $U'$ since $eU'\not\subseteq U$.
Moreover,
\[
 \dim_FU_1+\dim_FU_1'=\dim_FU+\dim_FU'
\]
by the dimension formula for a sum of subspaces. Also $\spn{U_1U_1'}\subseteq\spn{UU'}$.
Indeed, if $w\in U_1'$ then $Uw\subseteq\spn{UU'}$ and $eU'w=U'(ew)\subseteq\spn{UU'}$, since $ew\in U$.
Thus the claim for the pair $U,U'$ follows from the claim for the pair $U_1,U_1'$, which holds by the induction hypothesis.
Iterating the first assertion proves the bound for $\spn{U^j}$.
\end{proof}

By a derivation of a field $F$ we mean an additive map $\partial:F\to F$ such that $\partial(uw)=u\,\partial w+w\,\partial u$ for all $u,w\in F$.
We need derivations that detect independent directions among a specified set of generators. We give the field-theoretic argument here.

\begin{lemma}\label{lem:derivations}
Let $b$ be a perfect field of characteristic $p$, and let $F=b(z_1,\dots,z_N)$ be a finitely generated extension of $b$ of transcendence degree $h$.
Then there exist indices $i_1<\dots<i_h$ and derivations $\partial_1,\dots,\partial_h$ of $F$ that vanish on $b$ such that $\partial_lz_{i_l}=1$, and $\partial_lz_{i_m}=0$ for $m\ne l$.
\end{lemma}

\begin{proof}
Let $F^p:=\{u^p:u\in F\}$, which is a subfield of $F$ that contains $b$, since $b$ is perfect.
We first show that $[F:F^p]=p^h$.
Let $t_1,\dots,t_h$ be a transcendence basis of $F$ over $b$ and let $F_*:=b(t_1,\dots,t_h)$.
Then $F_*^p=b(t_1^p,\dots,t_h^p)$, and thus $[F_*:F_*^p]=p^h$.
Since the map $u\mapsto u^p$ is an isomorphism of $F$ onto $F^p$ that maps $F_*$ onto $F_*^p$, we have $[F^p:F_*^p]=[F:F_*]<\infty$.
Thus
\[
 [F:F^p]\,[F^p:F_*^p]=[F:F_*^p]=[F:F_*]\,[F_*:F_*^p]
\]
gives $[F:F^p]=p^h$.

Next, let $F_i:=F^p(z_1,\dots,z_i)$ for $0\le i\le N$.
Then $F_0=F^p$ and $F_N=F$.
The element $z_i$ is a root of the polynomial $T^p-z_i^p=(T-z_i)^p$, which has coefficients in $F^p$.
Thus its minimal polynomial over $F_{i-1}$ is $(T-z_i)^d$ for some $1\le d\le p$.
If $d<p$ then the coefficient $-dz_i$ of $T^{d-1}$ shows that $z_i\in F_{i-1}$.
It follows that $[F_i:F_{i-1}]$ is either $1$ or $p$, and it is equal to $p$ for exactly $h$ indices, say $i_1<\dots<i_h$.
Moreover, $F=F^p(z_{i_1},\dots,z_{i_h})$.

Now fix $l$ and let $F'$ be the field generated over $F^p$ by the elements $z_{i_m}$ with $m\ne l$.
Then $[F':F^p]\le p^{h-1}$ and $F=F'(z_{i_l})$.
It follows that $[F:F']=p$ and that the minimal polynomial of $z_{i_l}$ over $F'$ is $T^p-c$ with $c:=z_{i_l}^p$.
Thus $F$ is isomorphic to $F'[T]/(T^p-c)$.
The formal derivative $g\mapsto g'$ on $F'[T]$ maps the ideal generated by $T^p-c$ into itself, since $((T^p-c)g)'=(T^p-c)g'$.
It therefore induces a derivation $\partial_l$ of $F$ that vanishes on $F'$ and satisfies $\partial_lz_{i_l}=1$.
Since $F'$ contains $b$ and the elements $z_{i_m}$ with $m\ne l$, this completes the proof.
\end{proof}

We can now state the main result of this section.

\begin{proposition}\label{prop:dimension}
Let $b$ be a perfect field that contains $\F_p$, and let $F=b(z_0,\dots,z_{n-1})$ be a finitely generated extension of $b$ of transcendence degree $h\ge1$ such that $f$ is irreducible over $F$.
Suppose that $z=(z_0,\dots,z_{n-1})$ satisfies $\G(X(z))\in V_F$.
Then
\begin{equation}\label{eq:dimsum}
 \sum_{j=1}^k\min\{n,1+j(h-1)\}\le r.
\end{equation}
In particular, if $r<(k+1)n/2$ then $h\le r/D+1$.
\end{proposition}

\begin{proof}
Let $E:=K_F$, which is a field extension of $F$ of degree $n$, and let $X:=X(z)\in E$.
Let $\partial_1,\dots,\partial_h$ be the derivations given by Lemma~\ref{lem:derivations}.
For $c=(c_1,\dots,c_h)\in F^h$ the map $\partial_c:=c_1\partial_1+\dots+c_h\partial_h$ is a derivation of $F$ that vanishes on $b$.
We extend each derivation $\partial$ of $F$ that vanishes on $b$ to $E$ by letting it act on the coordinates with respect to the basis $1,\theta,\dots,\theta^{n-1}$, and to $E^k$ by letting it act on each component.
The extension is a derivation of $E$, since the products $\theta^i\theta^{i'}$ are linear combinations of the elements of this basis with coefficients in $\F_p$, on which $\partial$ vanishes.
Let
\[
 U:=\{\partial_cX:c\in F^h\}\subseteq E.
\]
This is an $F$-subspace of $E$ of dimension $h$. Indeed, the coefficient of $\theta^{i_m}$ in $\partial_c X$ is $\partial_c z_{i_m}=c_m$, so $c\mapsto\partial_c X$ is injective.

Let $\partial=\partial_c$ for some $c\in F^h$.
We claim that
\begin{equation}\label{eq:derivinV}
 \partial^j\G(X)\in FV_0\qquad\text{for every }j\ge1.
\end{equation}
Indeed, let $\ell:\F_q^k\to\F_p$ be a linear functional that vanishes on $V_0$.
In the basis of $E^k$ over $F$ formed by the vectors with one coordinate equal to a power of $\theta$ and the others equal to zero, $\ell_F$ is a linear form with coefficients in $\F_p$.
Thus $\ell_F(\partial w)=\partial(\ell_F(w))$ for every $w\in E^k$.
Since $\ell_F(\G(X))=\ell(v)$ belongs to $\F_p$ by \eqref{eq:functionals}, we obtain $\ell_F(\partial^j\G(X))=\partial^j(\ell(v))=0$ for $j\ge1$, and \eqref{eq:derivinV} follows from \eqref{eq:functionals}, applied to $V_0$ in place of $V$.

We now define a map $\Phi:E^k\to E^k$ by
\[
 (\Phi w)_i:=\sum_{l=1}^i\binom il(-X)^{i-l}w_l\qquad(1\le i\le k).
\]
This map is linear over $E$ and is given by a triangular matrix with ones on the diagonal.
Thus it is a bijection, and $W:=\Phi(FV_0)$ is a subspace of $E^k$ over $F$ with $\dim_FW=r$.
We claim that for $1\le j\le k$ the vector $w^{(j)}:=\Phi(\partial^j\G(X))$, which belongs to $W$ by \eqref{eq:derivinV}, satisfies
\begin{equation}\label{eq:triangular}
 w^{(j)}_j=j!\,(\partial X)^j\qquad\text{and}\qquad w^{(j)}_i=0\quad\text{for }j<i\le k.
\end{equation}
To see this, we extend $\partial$ to a derivation of the polynomial ring $E[Z]$ by letting it act on the coefficients.
By the binomial theorem, and since $\partial^j1=0$, we have
\[
 \partial^j\bigl((X-Z)^i\bigr)=\sum_{l=1}^i\binom il(-Z)^{i-l}\partial^j(X^l),
\]
and the value of the right-hand side at $Z=X$ is $w^{(j)}_i$.
On the other hand, by the product rule $\partial^j((X-Z)^i)$ is a sum of terms of the form
\[
 \frac{j!}{j_1!\cdots j_i!}\prod_{m=1}^i\partial^{j_m}(X-Z)
\]
with $j_1+\dots+j_i=j$.
A term with $j_m=0$ for some $m$ contains the factor $X-Z$, and vanishes at $Z=X$.
If $i>j$ then every term is of this type.
If $i=j$ the only other term is the one with $j_1=\dots=j_i=1$, which is equal to $j!(\partial X)^j$.
This proves \eqref{eq:triangular}.

For $0\le j\le k$ let $W_j$ denote the space of all $w\in W$ with $w_i=0$ for $i>j$, so that $\{0\}=W_0\subseteq W_1\subseteq\dots\subseteq W_k=W$.
The map $w\mapsto w_j$ from $W_j$ to $E$ has kernel $W_{j-1}$.
Denoting its image by $P_j$, we thus have
\begin{equation}\label{eq:filtration}
 r=\dim_FW=\sum_{j=1}^k\dim_FP_j.
\end{equation}
By \eqref{eq:triangular} the space $P_j$ contains $j!u^j$ for every $u\in U$.
Since $p>k$, the integer $j!$ is invertible in $\F_p$, and the identity
\[
 j!\,u_1u_2\cdots u_j=\sum_{I\subseteq\{1,\dots,j\}}(-1)^{j-|I|}\Bigl(\sum_{i\in I}u_i\Bigr)^j
\]
shows that $P_j$ contains the product of any $j$ elements of $U$.
Thus $P_j\supseteq\spn{U^j}$, and Lemma~\ref{lem:products} gives $\dim_FP_j\ge\min\{n,1+j(h-1)\}$.
Together with \eqref{eq:filtration} this proves \eqref{eq:dimsum}.

Finally, suppose that $r<(k+1)n/2$ and $h-1>r/D$.
Then $1+j(h-1)>jr/D$, and $n>kr/D\ge jr/D$ for $1\le j\le k$.
Thus the left-hand side of \eqref{eq:dimsum} exceeds $\sum_{j=1}^kjr/D=r$, which is a contradiction.
This completes the proof.
\end{proof}

\section{Counting points}\label{sec:counting}

We keep the notation of the previous section.
For $m\ge1$ let $N_m$ denote the number of solutions $z\in\F_{p^m}^n$ of \eqref{eq:condition} with $F=\F_{p^m}$.
In other words, $N_m$ is the number of solutions over $\F_{p^m}$ of the system of polynomial equations described after \eqref{eq:condition}.
Our goal is to estimate $N_1$.
We let $\Fbar$ denote an algebraic closure of $\F_p$, and we regard the fields $\F_{p^m}$ as subfields of $\Fbar$.
Let
\[
 b:=\bigcup_{n\nmid m}\F_{p^m}\subseteq\Fbar.
\]

\begin{lemma}\label{lem:fieldb}
The set $b$ is a perfect field, and $f$ is irreducible over $b$.
\end{lemma}

\begin{proof}
If $n$ divides neither $m$ nor $m'$ then it does not divide their least common multiple $m''$, since $n$ is prime, and $\F_{p^m}$ and $\F_{p^{m'}}$ are contained in $\F_{p^{m''}}$.
Thus $b$ is a field.
It is perfect since each of the fields $\F_{p^m}$ is.
Let $n\nmid m$.
The roots of $f$ in $\Fbar$ are $\theta^{p^i}$, $0\le i<n$, and the map $u\mapsto u^{p^m}$ permutes them cyclically with a shift by $m$.
Since $n$ is prime and does not divide $m$ this permutation is transitive, and thus $f$ is irreducible over $\F_{p^m}$.
A factorization of $f$ over $b$ would be a factorization over one of these fields.
\end{proof}

Let $S\subseteq b^n$ denote the set of solutions $z\in b^n$ of \eqref{eq:condition} with $F=b$.
For $n\nmid m$ we have $N_m=|S\cap\F_{p^m}^n|$.
We will bound this number by means of an elementary form of dimension theory for the set $S$, which we now describe.
For a set $P\subseteq b^n$ we let $I(P)$ denote the ideal of all polynomials in $b[T_0,\dots,T_{n-1}]$ that vanish on $P$.
We say that $P$ is irreducible if it is nonempty and $I(P)$ is a prime ideal.
In this case $b[P]:=b[T_0,\dots,T_{n-1}]/I(P)$ is an integral domain, whose elements we regard as functions on $P$ with values in $b$.
We let $F_P$ denote its field of fractions.
Thus $F_P=b(z_0,\dots,z_{n-1})$, where $z_i$ denotes the image of $T_i$, that is, the $i$-th coordinate function on $P$.
The dimension of $P$ is defined to be the transcendence degree of $F_P$ over $b$.

\begin{lemma}\label{lem:decomposition}
Every subset of $b^n$ is the union of finitely many irreducible sets.
\end{lemma}

\begin{proof}
Here the empty set is regarded as the union of an empty collection of sets.
Suppose that the lemma fails.
By Hilbert's basis theorem (see \cite[Theorem~7.5]{AM}) every nonempty collection of ideals of $b[T_0,\dots,T_{n-1}]$ has a maximal element.
Thus there is a set $P\subseteq b^n$ that is not a finite union of irreducible sets and for which $I(P)$ is maximal among the ideals of such sets.
The set $P$ is nonempty and is not irreducible.
Thus there are polynomials $u$ and $w$ that do not belong to $I(P)$ with $uw\in I(P)$.
Then $P$ is the union of the sets $P\cap\{u=0\}$ and $P\cap\{w=0\}$.
The ideals of these two sets contain $u$ and $w$, respectively, and thus they are strictly larger than $I(P)$.
It follows that both sets are finite unions of irreducible sets, so the same holds for $P$.
This is a contradiction.
\end{proof}

\begin{lemma}\label{lem:closed}
Let $P\subseteq b^n$ be an irreducible set.
Then every element of $F_P$ that is algebraic over $b$ belongs to $b$, and $f$ is irreducible over $F_P$.
\end{lemma}

\begin{proof}
Let $\phi\in F_P$ be algebraic over $b$, with minimal polynomial $g\in b[Y]$ of degree $d$.
We write $\phi=u/w$ with $u,w\in b[P]$ and $w\ne0$.
Since $w$ is a nonzero function on $P$, there is a point $a\in P$ with $w(a)\ne0$.
The element $w^dg(u/w)$ of $b[P]$ vanishes, and evaluating it at $a$ we find that $g(u(a)/w(a))=0$.
Thus the irreducible polynomial $g$ has a root in $b$, which gives $d=1$ and $\phi\in b$.

Now suppose that $f=f_1f_2$, where $f_1$ and $f_2$ are monic polynomials with coefficients in $F_P$.
The coefficients of $f_1$ are polynomial expressions in the roots of $f$ in an extension of $F_P$, and these roots are algebraic over $\F_p$.
Thus the coefficients of $f_1$ are algebraic over $b$, and by the first part of the proof they belong to $b$.
Since $f$ is irreducible over $b$ by Lemma~\ref{lem:fieldb}, one of $f_1$ and $f_2$ is constant.
\end{proof}

We can now apply Proposition~\ref{prop:dimension} to the coordinate functions on each irreducible set.

\begin{corollary}\label{cor:pieces}
Suppose that $r<(k+1)n/2$.
Then $S$ is the union of finitely many irreducible sets of dimension at most $r/D+1$.
\end{corollary}

\begin{proof}
By Lemma~\ref{lem:decomposition} it suffices to show that every irreducible set $P\subseteq S$ has dimension at most $r/D+1$.
Let $F:=F_P$ and let $z=(z_0,\dots,z_{n-1})\in F^n$ be the vector of the coordinate functions on $P$.
The polynomials in the system that expresses the condition \eqref{eq:condition} vanish on $S$, and thus on $P$.
Thus they vanish at $z$, which means that $\G(X(z))\in V_F$.
The field $b$ is perfect by Lemma~\ref{lem:fieldb}, and $f$ is irreducible over $F$ by Lemma~\ref{lem:closed}.
Thus, if the transcendence degree $h$ of $F$ over $b$ is positive, Proposition~\ref{prop:dimension} gives $h\le r/D+1$.
\end{proof}

We next count the points whose coordinates lie in a finite subfield of $b$. The constant in this count may depend on the equations defining $S$.

\begin{lemma}\label{lem:box}
Let $S\subseteq b^n$ be the union of finitely many irreducible sets of dimension at most $H$, where $H\ge0$ is an integer.
Then there is a constant $C_S$ such that $|S\cap\F_Q^n|\le C_SQ^H$ for every finite subfield $\F_Q$ of $b$.
\end{lemma}

\begin{proof}
If $H\ge n$, the estimate follows from $|\F_Q^n|=Q^n$. We may therefore assume $H<n$.
Let $I\subseteq\{0,\dots,n-1\}$ be a set of $H+1$ indices.
On each of the irreducible sets $P$ in the statement the coordinate functions $z_i$, $i\in I$, are algebraically dependent over $b$.
Thus there is a nonzero polynomial in the variables $T_i$, $i\in I$, with coefficients in $b$ that vanishes on $P$.
Let $g_I$ denote the product of these polynomials over the sets $P$.
Then $g_I$ is a nonzero polynomial in the variables $T_i$, $i\in I$, that vanishes on $S$.
Let $B$ be an integer that exceeds the degrees of all of the polynomials $g_I$.

We fix a monomial order, for instance the lexicographic order.
Recall that this is a well-ordering of the set of monomials in $T_0,\dots,T_{n-1}$ that is compatible with multiplication.
Let $S_Q:=S\cap\F_Q^n$.
The functions from $S_Q$ to $b$ form a vector space over $b$ of dimension $|S_Q|$, and this space is spanned by the restrictions of the monomials to $S_Q$.
Indeed, the function that is equal to one at a point $a\in\F_Q^n$ and vanishes at the other points of $\F_Q^n$ is the restriction of the polynomial $\prod_i(1-(T_i-a_i)^{Q-1})$.
We say that a monomial is reduced if all of its exponents are less than $Q$ and at most $H$ of them are at least $B$.
We claim that the restriction of every monomial to $S_Q$ is a linear combination of restrictions of reduced monomials.
Suppose that this fails, and consider the smallest monomial $T^\beta=T_0^{\beta_0}\cdots T_{n-1}^{\beta_{n-1}}$ for which it fails.
If $\beta_i\ge Q$ for some $i$, then $T^\beta$ agrees on $\F_Q^n$ with the smaller monomial $T^\beta T_i^{-(Q-1)}$.
The claim holds for the latter, so it also holds for $T^\beta$, which is a contradiction.
Otherwise, since $T^\beta$ is not reduced, there is a set $I$ of $H+1$ indices such that $\beta_i\ge B$ for all $i\in I$.
Let $T^\delta$ denote the largest monomial that occurs in $g_I$.
Then $T^\delta$ divides $T^\beta$, by the choice of $B$.
Since $g_I$ vanishes on $S$, the restriction of $T^\delta$ to $S$ is a linear combination of restrictions of smaller monomials.
Multiplying by $T^{\beta-\delta}$ we see that the restriction of $T^\beta$ to $S_Q$ is a linear combination of restrictions of monomials that are smaller than $T^\beta$.
The claim holds for each of these, so it also holds for $T^\beta$, which is again a contradiction.
This proves the claim.

It follows that $|S_Q|$ is at most the number of reduced monomials.
A reduced monomial is determined by a set of $H$ indices that contains all indices $i$ with $\beta_i\ge B$, the exponents at these indices, and the remaining exponents, which are less than $B$.
Thus $|S_Q|\le\binom nHQ^HB^{n-H}$.
\end{proof}

The constant $C_S$ in Lemma~\ref{lem:box} depends on the set $S$ in a way that we do not control.
In order to remove it we use the zeta function.
For the moment let $g_1,\dots,g_M$ be arbitrary polynomials in $N$ variables with coefficients in $\F_p$, and let $N_m$ denote the number of their common zeros in $\F_{p^m}^N$.
By a theorem of Dwork \cite{Dwork} the power series $\exp(\sum_{m\ge1}N_mt^m/m)$ is a rational function of $t$.
Writing it in lowest terms as a quotient of polynomials with constant term one, and comparing logarithmic derivatives, one finds that there exist nonzero complex numbers $\omega_1,\dots,\omega_\Lambda$ and signs $\epsilon_1,\dots,\epsilon_\Lambda\in\{1,-1\}$ such that
\begin{equation}\label{eq:dwork}
 N_m=\sum_{i=1}^\Lambda\epsilon_i\omega_i^m\qquad\text{for all }m\ge1.
\end{equation}
Here $\Lambda$ is the sum of the degrees of the numerator and the denominator of the zeta function.
The following bound for this number is due to Wan and Zhang \cite{WZ}.
It is important for us that it does not depend on the number $M$ of equations.

\begin{theorem}[Wan and Zhang]\label{thm:WZ}
Suppose that the polynomials $g_1,\dots,g_M$ have degree at most $d$.
Then the number $\Lambda$ in \eqref{eq:dwork} satisfies
\[
 \Lambda\le3^{N+1}\binom{2N}N(2d+1)^N<3(24d+12)^N.
\]
\end{theorem}

More precisely, \cite[Corollary~1.1.4]{WZ} states that the right-hand side bounds the sum of the compactly supported $\ell$-adic Betti numbers of the variety defined by $g_1,\dots,g_M$ over $\Fbar$, and $\Lambda$ is at most this sum (see \cite[Section~1.1]{WZ}).
The earlier bounds of Bombieri \cite{Bombieri} and Katz \cite{Katz} depend also on the number of equations. In our application there are at most $kn$ equations, so either bound would give Theorem~\ref{thm:section} with a larger constant in place of $24k+12$.

\begin{lemma}\label{lem:zeta}
Let $N_m$ and $\Lambda$ be as in \eqref{eq:dwork}.
Let $n\ge1$ and $H\ge0$ be integers, and suppose that there is a constant $C$ such that $N_m\le Cp^{mH}$ for all $m\equiv1\pmod n$.
Then $N_1\le\Lambda p^H$.
\end{lemma}

\begin{proof}
For each complex number $\beta$ let $c_\beta$ denote the sum of $\epsilon_i\omega_i$ over the indices $i$ with $\omega_i^n=\beta$.
By \eqref{eq:dwork} we have $N_{1+nt}=\sum_\beta c_\beta\beta^t$ for every integer $t\ge0$, and thus
\[
 \sum_{t\ge0}N_{1+nt}\,y^t=\sum_\beta\frac{c_\beta}{1-\beta y}
\]
for small $y$.
By hypothesis the power series on the left converges for $|y|<p^{-nH}$.
The right-hand side is a rational function with a pole at $y=1/\beta$ for each $\beta$ with $c_\beta\ne0$.
Thus $|\beta|\le p^{nH}$ whenever $c_\beta\ne0$.
For such $\beta$ every $\omega_i$ with $\omega_i^n=\beta$ satisfies $|\omega_i|\le p^H$, which gives $|c_\beta|\le p^H\,|\{i:\omega_i^n=\beta\}|$.
Taking $t=0$ and summing these bounds gives $N_1=|\sum_\beta c_\beta|\le\Lambda p^H$.
\end{proof}

\begin{proof}[Proof of Theorem~\ref{thm:section}]
Let $H$ be the largest integer with $H\le r/D+1$.
By Corollary~\ref{cor:pieces} and Lemma~\ref{lem:box} there is a constant $C$ such that $N_m\le Cp^{mH}$ for every $m$ that is not divisible by $n$, and in particular for every $m\equiv1\pmod n$.
The numbers $N_m$ count the solutions over $\F_{p^m}$ of a system of polynomial equations of degree at most $k$ in $n$ variables with coefficients in $\F_p$.
Thus Lemma~\ref{lem:zeta} and Theorem~\ref{thm:WZ} give
\[
 N_1\le\Lambda p^H<3(24k+12)^n\,p^{r/D+1}.
\]
Since $N_1$ is the number of $x\in\F_q$ with $\G(x)\in V$, this completes the proof.
\end{proof}

\section{Lower bounds for the extractor error}\label{sec:lower}

For a map $E:\Omega\times\Omega\to\Omega'$ let $\varepsilon_E(N)$ be the largest total variation distance from uniform when the two inputs are independent and each is uniform on a set of exactly $N$ points. Every probability measure with atoms of mass at most $1/N$ is a convex combination of uniform $N$-point measures. Thus $\varepsilon_E(N)$ is also the worst error at min-entropy $\log_2 N$. For an output of the form $f(xy)$ on a field, a constant output on $A\times B$ persists on every $tA\times t^{-1}B$. This lets us use many copies of the same example when $N$ is larger than $|A|$ and $|B|$.

\subsection{Fixed characteristic and growing extension degree}\label{sec:lowerfixed}

Fix a prime $p>k$ and let $q=p^n$. Here $n$ may tend to infinity through arbitrary positive integers. We consider a linear projection of $P_k(xy)$, as in Corollary~\ref{cor:extractor}.

\begin{proposition}\label{prop:fixedoutput}
Let $D=k(k+1)/2$ and let $\pi:\F_q\to\F_p^m$ be any surjective linear map with $m\ge1$. For $E(x,y)=\pi(P_k(xy))$ and $1\le N\le q/2$,
\begin{equation}\label{eq:fixedlower}
  \varepsilon_E(N)\gtrsim_{p,k}\min\{1,q^{1/(D+1)}/N\}.
\end{equation}
For the coordinate giving the coefficient of $\theta^{n-1}$ in a power basis, the numerator can be replaced by $q^{1/(2k)}$.
\end{proposition}

\begin{proof}
Take any nonzero linear functional on the output and write $\ell$ for its composition with $\pi$. For $y\in I_h$, the function $y\mapsto\ell(P_k(xy))$ is determined by
\[
 \ell(x^j\theta^i),\qquad 1\le j\le k,\quad 0\le i\le j(h-1).
\]
Indeed $y^j$ is a polynomial in $\theta$ of degree at most $j(h-1)$. The list has $D(h-1)+k$ entries. Some list occurs for at least $q p^{-D(h-1)-k}$ values of $x$. On this set the output coordinate depends only on $y$, so it is constant on some subset of $I_h$ of size at least $p^{h-1}$. With $h=1+\lfloor n/(D+1)\rfloor$ both sets have size $\gtrsim_{p,k}q^{1/(D+1)}$.

For the specified power-basis coordinate use $I_h$ for both sets, with $h=1+\lfloor(n-2)/(2k)\rfloor$. The ordinary polynomial representing $P_k(xy)$ has degree at most $2k(h-1)\le n-2$, so its last coefficient is zero. Each set has size $\gtrsim_p q^{1/(2k)}$.

In both cases we may discard zero. If the coordinate is constant on $A\times B$, it has the same value on $tA\times t^{-1}B$ for every $t\ne0$. Lemma~\ref{lem:copies} shows that sets of size $M$ with this property give error $\gtrsim\min\{1,M/N\}$ at size $N$. This proves the proposition.
\end{proof}

For $k=2$ both exponents are $1/4$. At $N=\lceil q^\rho\rceil$ the error of every nonzero linear coordinate is therefore $\gtrsim_p q^{1/4-\rho}$ for $1/4\le\rho<1$. It stays bounded away from zero at rate $1/4$ and is at least a constant times $q^{-1/12}$ at rate $1/3$. The sources may depend on the coordinate. Proposition~\ref{prop:affineoutput} gives the exact threshold if both sources must be affine subspaces and the coordinate may be chosen.

\subsection{Prime fields as \texorpdfstring{$p\to\infty$}{p tends to infinity}}\label{sec:lowerprime}

For $q=p$, a nonzero linear coordinate retains the whole field value, up to a nonzero scalar. Already $P_2(t)=(t+1/2)^2-1/4$ takes only $(p+1)/2$ values. Thus its field-valued output has error at least $(p-1)/(2p)$ even on uniform inputs. Now hold $n$ fixed and let $p$ grow. Two independent inputs uniform on the subfield $\F_p\subseteq\F_{p^n}$ each have rate $1/n$. On this subfield $\ell(z)=z\ell(1)$ for every $\F_p$-linear coordinate $\ell$. If $\ell(1)\ne0$, the same lower bound applies. If $\ell(1)=0$, the output is identically zero and its error is exactly $(p-1)/p$.

The prime-field constructions of Bourgain and of Chattopadhyay and Gurumukhani take bits of the least nonnegative residue $[P_k(xy)]_p$. We consider the half-interval bit $1_{\{[z]_p<p/2\}}$ and the parity bit $[z]_p\bmod 2$. A lower bound for the parity bit also applies to any output that includes it.

\begin{proposition}\label{prop:primeoutput}
Fix $k\ge2$. Let $p$ be a sufficiently large prime, let $\Omega=\{0,\dots,Q-1\}$ with $p/2\le Q\le p$, and let $E$ be either bit of $P_k(xy)$ on $\Omega\times\Omega$. For $1\le N\le Q/2$,
\[
 \varepsilon_E(N)\gtrsim_k\min\{1,p^{1/(2k)}/N\}.
\]
For either bit of $\lambda P_k(xy)$, uniformly in $\lambda\in\F_p^\times$, the same estimate holds with $p^{1/(D+1)}$ in the numerator.
\end{proposition}

\begin{proof}
Take $A=B$ to be the positive even integers up to $c_k p^{1/(2k)}$, with $c_k$ sufficiently small. For $x,y\in A$ the integer $P_k(xy)$ is even and belongs to $[0,p/2)$. This is the interval example in \cite[Section~5]{ChG}. Lemma~\ref{lem:copies} gives the stated bound for $\Omega=\F_p$. For the truncated input interval, $P_k(-xy)$ is also even, has sign $(-1)^k$, and has absolute value less than $p/2$. When the sign is negative its residue is $p+P_k(-xy)$, so each bit is again constant. The last part of the lemma then applies. The construction for a general multiplier is given in Lemma~\ref{lem:multiplier}.
\end{proof}

On the full prime field, take $N=\lceil p^\rho\rceil$ with $1/4\le\rho<1$. For $xy+x^2y^2$ both bits have worst error $\gtrsim p^{1/4-\rho}$. Restricting the inputs to an initial interval of length $Q\ge p/2$ changes the input length by at most one bit. An error guarantee $p^{-\delta}$ with fixed $\delta>0$ therefore requires $\rho\ge1/4+\delta-O(1/\log p)$. For the scalar quadratic the range $1/4<\rho\le1/3$ remains open. Over extensions we ask this for each fixed $\rho$, all sufficiently large fixed $p$ and prime degrees $n\to\infty$.

For comparison, the paraboloid in $\F_p^3$ gives the field value $x\cdot y+(x\cdot x)(y\cdot y)$ for $x,y\in\F_p^2$. Its bit can be constant on two sets of size $\gtrsim p^{2/3}$, giving constant error at rate $1/3$. Appendix~\ref{app:paraboloids} gives the proof, the explicit single-set example of size $\asymp\sqrt p\log p$ for $p\equiv3\pmod8$, and the $5/12$ example for inputs in $\F_p^3$.

\section{The moment order and the characteristic}\label{sec:remarks}

Can one take $s=D=k(k+1)/2$ in the moment estimate of Theorem~\ref{thm:main}? Here $s$ counts the summands on each side of the $k$ equations in \eqref{eq:system}. For every $1/(k+1)<\rho\le1$ and every prime $p$ satisfying \eqref{eq:plarge}, we ask whether there is an $\eta=\eta(k,\rho,p)>0$ for which
\[
 J_{D,k}(A)\le |A|^{2D}q^{-k/2-\eta}
\]
for every $A\subseteq\F_q$ with $|A|\ge q^\rho$, whenever $n$ is a sufficiently large prime. The diagonal solutions give $J_{s,k}(A)\ge|A|^s$. Applying the desired estimate to a set of size $\lceil q^\rho\rceil$ forces $s\rho\ge k/2+\eta$. Thus $D$ is the smallest integer $s$ that could work uniformly as $\rho$ decreases to $1/(k+1)$. For the quadratic system this asks for three summands on each side, or six field variables in all.

The estimate with $s=D$ does hold for subsets of certain affine images of $I_h$.
Let $h\ge1$ be such that $k(h-1)<n$ and let $A\subseteq I_h$.
We identify $A$ with a set $A'$ of polynomials of degree less than $h$ in $\F_p[T]$.
If $x_1,\dots,x_s,y_1,\dots,y_s\in A'$ then the two sides of the $j$-th equation in \eqref{eq:system} are polynomials of degree less than $n$, and they agree at $\theta$ only if they are equal.
Thus $J_{s,k}(A)=J_{s,k}(A')$.
We use Wooley's function-field mean value theorem \cite[Theorem~17.1]{WooleyNEC}. In Section~17 of that paper he sketches the adaptation of the efficient congruencing argument to function fields and refers to the integer case for the details of that step. Taking the coefficients to be the characteristic function of $A'$ gives
\[
 J_{D,k}(A')\lesssim_{p,k,\eps}p^{\eps h}|A'|^D
\]
for every $\eps>0$, provided that $p>k$.
The count $J_{s,k}(A)$ is invariant under the maps $x\mapsto ax+a'$ with $a\ne0$. Indeed the binomial theorem expresses each translated equation as a linear combination of equations of the same or lower degree, and the affine map is invertible.
It follows that if $A$ is contained in a set of the form $aI_h+a'$ with $k(h-1)<n$, and $|A|\ge q^\rho$ for some $\rho>\frac1{k+1}$, then
\[
 J_{D,k}(A)\lesssim_{p,k,\eps}|A|^{2D}q^{\eps-D\rho}.
\]
Since $q^\rho\le|A|\le p^h$ and $k(h-1)<n$, these hypotheses require $\rho\le h/n<1/k+1/n$. Choose $0<\eps<D\rho-k/2$ and then $0<\eta<D\rho-k/2-\eps$. This gives \eqref{eq:goal} with $s=D$ for all sufficiently large $q$ in each fixed characteristic $p>k$, with no restriction on the degree $n$.
When $D$ is odd, monotonicity gives the same normalized estimate with $s=D+1$ for use in \eqref{eq:bilinear}.
For unrestricted sets $A$ we do not know even whether some $s=s(k)$ works for all the admissible $\rho,p$.

A related question concerns the characteristic. The lower bound \eqref{eq:plarge} on $p$ enters only through the factor $(24k+12)^n$ in Theorem~\ref{thm:section}.
Replacing this factor by $q^{o(1)}$ for fixed $k,p$ would give the moment estimate for every $p>k$ and $\rho>1/(k+1)$.
For $p\le k$ there are additional relations. The equation of degree $pj$ is the $p$-th power of the equation of degree $j$, so only $k-\lfloor k/p\rfloor$ of the equations need to be imposed.
In characteristic two and even degree $k$, this alone gives $J_{s,k}(A)\ge |A|^{2s}q^{-k/2}$ for every nonempty $A$. For the quadratic extractor in characteristic two, $\Tr(t+t^2)=0$ for every $t$.

The proof of Theorem~\ref{thm:section} also gives a bound for affine subspaces of dimension $r\ge(k+1)n/2$. For any $r$, the number of $x$ with $\G(x)\in V$ is at most $3(24k+12)^np^H$, where $H$ is the largest of $0$ and the integers $1\le h\le n$ satisfying \eqref{eq:dimsum}.

Over prime fields the corresponding question asks for uniformity in the characteristic. Given $k\ge2$ and $1/(k+1)<\rho\le1$, do there exist $s=s(k,\rho)$ and $\eta=\eta(k,\rho)>0$ such that
\[
 J_{s,k}(A)\le |A|^{2s}p^{-k/2-\eta}
\]
for all sufficiently large primes $p$ and all $A\subseteq\F_p$ with $|A|\ge p^\rho$?
Intervals give the necessary threshold, and \eqref{eq:quadratic} answers the question for $k=2$. We do not know how to obtain the full range for any $k\ge3$.
The dependence on $p$ in Theorem~\ref{thm:Marton} is harmless when we fix $p$ and let $n$ grow. When $p$ itself grows, the covering bound used in Proposition~\ref{prop:inverse} is no longer uniform. Bourgain's prime-field argument uses sum-product estimates at this point.

\appendix
\section{Proofs for the scalar lower bounds}\label{app:lowerproofs}

We first record two ways to pass from sets on which the output is restricted to flat sources of a prescribed size. In the following statements, write $U_C$ for the uniform measure on a nonempty finite set $C$.

\begin{lemma}\label{lem:enlarge}
Suppose $E(A\times B)\subseteq T\subseteq\Omega'$ for nonempty subsets of $\Omega$. For $N\le|\Omega|/2$,
\[
 \varepsilon_E(N)\ge\frac{1-|T|/|\Omega'|}{4}
       \min\{1,|A|/N\}\min\{1,|B|/N\}.
\]
\end{lemma}

\begin{proof}
For $C\subseteq\Omega$ with $|C|\le N$ choose a set $S$ of $N$ points outside $C$. Let $\mu$ be the average of the uniform measures on $C\cup S'$, where $S'$ ranges over the $(N-|C|)$-element subsets of $S$. Then
\[
 \frac{1_C}{N}=\mu-(1-|C|/N)U_S.
\]
Consequently, for $g_T(x,y)=1_T(E(x,y))-|T|/|\Omega'|$ and sets $C,C'$ of size at most $N$,
\begin{equation}\label{eq:smallrect}
 \Big|\sum_{x\in C,\,y\in C'}g_T(x,y)\Big|\le 4N^2\varepsilon_E(N).
\end{equation}
Indeed each of the four pairs of probability measures on the right after expansion is a convex combination of pairs of uniform $N$-point measures. Apply \eqref{eq:smallrect} to $A$ and $B$, truncating either set that has more than $N$ points.
\end{proof}

For an output depending only on $xy$, copies of the same example give a stronger estimate. The two inputs in the resulting example are independent.

\begin{lemma}\label{lem:copies}
Let $K$ be a finite field of odd order and let $E(x,y)=f(xy)$ on $K\times K$ with output in $\Omega'$. Suppose $A,B\subseteq K^\times$ have size $M$ and $f(AB)\subseteq T$, where $|T|\le|\Omega'|/2$ and $AB=\{ab:a\in A,b\in B\}$. Then for $1\le N\le |K|/2$,
\[
 \varepsilon_E(N)\gtrsim\min\{1,M/N\}.
\]
If $K=\F_p$ and $f(-AB)$ is also contained in a set of at most $|\Omega'|/2$ outputs, this conclusion holds on $\Omega=\{0,\dots,Q-1\}$ for $p/2\le Q\le p$ and $N\le Q/2$, provided $M\le Q/2$.
\end{lemma}

\begin{proof}
First suppose $M\le N$. Put $r=\max\{1,\lfloor N/(2M)\rfloor\}$, so $r\gtrsim N/M$ and $rM\le N$. We choose $r$ pairs of sets $A_i\subseteq t_iA$, $B_i\subseteq t_i^{-1}B$, with the $A_i$ pairwise disjoint and the $B_i$ pairwise disjoint. At each choice the previously used points on either side number at most $(r-1)M\le N/2\le |K^\times|/4$. For a uniform $t\in K^\times$ the expected total number of used points in $tA$ and $t^{-1}B$ is at most $M/2$. Some $t$ therefore leaves at least $M/2$ unused points on each side. Take these as $A_i,B_i$. On each $A_i\times B_i$ the output belongs to $T$.

Let $\sigma_1,\dots,\sigma_r$ be independent uniform signs. In the sum
\[
 \sum_{x,y}\Big(\sum_i\sigma_i1_{A_i}(x)\Big)
            \Big(\sum_j\sigma_j1_{B_j}(y)\Big)
            \big(1_T(E(x,y))-|T|/|\Omega'|\big)
\]
averaging over the signs retains just the pairs with $i=j$. The average is at least $\frac12\sum_i|A_i||B_i|\gtrsim rM^2$. For each choice of signs, each signed sum is the difference of the characteristic functions of two sets of size at most $N$. By \eqref{eq:smallrect} the absolute value of the displayed sum is at most $16N^2\varepsilon_E(N)$. This gives $\varepsilon_E(N)\gtrsim M/N$.

For the truncated interval, every nonzero residue has at least one of its two signs in $\Omega$. For each $i$, some choice of signs $a_i,b_i\in\{1,-1\}$ puts at least half of $A_i$ and half of $B_i$ in $\Omega$ after multiplying by $a_i$ and $b_i$, respectively. Keep those points and the most frequent pair of signs. There are at least $r/4$ pairs left, disjoint on each side because the two signs are now fixed. Their products all lie in $AB$ or all in $-AB$, so the same averaging argument applies in $\Omega$. Finally $\varepsilon_E(N)$ is nonincreasing in $N$, since a uniform measure on a larger set is an average of uniform measures on its $N$-point subsets. For $N<M$ one may truncate the two original sets directly in the full field. In the restricted interval apply monotonicity with the already proved case $N=M$.
\end{proof}

\begin{lemma}\label{lem:multiplier}
For either of the prime-field bits, every $\lambda\ne0$ and every sufficiently large prime $p$, there are $A,B\subseteq\F_p^\times$ of size $\gtrsim_k p^{1/(D+1)}$ such that the bit of $\lambda P_k(xy)$ and the bit of $\lambda P_k(-xy)$ are each constant on $A\times B$.
\end{lemma}

\begin{proof}
Put $H=\lfloor p^{1/(D+1)}\rfloor$ and $Y=\{2,4,\ldots,2H\}$. Fix a sufficiently small $\delta>0$ in terms of $k$. Partition $x\in\F_p^\times$ according to the intervals containing $[\lambda x^j]_p$, $1\le j\le k$, using consecutive intervals of length $\lfloor\delta p/(k(2H)^j)\rfloor$ and a possibly shorter last interval. These lengths grow with $p$ since $k<D+1$. There are $O_k(H^D)$ sets in the partition. For fixed $y\in Y$ and either sign the even integers
\[
                  \sum_{j=1}^k[\lambda x^j]_p(\pm y)^j
\]
vary by at most $\delta p$ within any one set. If one of their residues has distance greater than $\delta p$ from $0,p/2,p$, both its half and its parity are fixed on that set. Indeed, the integer quotient on division by the odd prime $p$ is unchanged there.

For each $y$, at least $3(p-1)/4$ values of $x$ satisfy this distance condition for both signs. For either sign each residue has at most $k$ preimages, and the excluded intervals altogether contain at most $6\delta p+6$ residues. Discard the partition sets of size less than $(p-1)/(4J)$, where $J$ is the number of sets. In the remaining sets there are still at least $(p-1)|Y|/2$ pairs satisfying both conditions. Thus one such set $A$ contains an admissible $x$ for at least $|Y|/2$ different $y$. For each such $y$ the two signs give a fixed pair of output bits across $A$. Keeping the most frequent pair gives $|B|\ge |Y|/8$. Since $|A|\gtrsim_k p/H^D$, both sizes have the required order.
\end{proof}

The same argument works for any finite group acting on the two input spaces and preserving $E$, provided each orbit meeting a source has at least $N$ points. In the proof of Lemma~\ref{lem:copies} take $r=\max\{1,\lfloor N/(4M)\rfloor\}$. A uniform group element sends each point uniformly onto its orbit. The probability of hitting a previously used point on its side is at most $(r-1)M/N\le1/4$, which is the only estimate on the action used in the proof.

When the two sources are restricted to affine subspaces, the least possible dimension of the output can be computed exactly.

\begin{proposition}\label{prop:affineoutput}
Let $p>k$, let $n$ be prime and let $q=p^n$. Among affine $\F_p$-subspaces $A,B\subseteq\F_q$ of dimensions $h_A,h_B\ge1$, the smallest possible dimension of the affine span of $\{P_k(xy):x\in A,\ y\in B\}$ is
\[
 \min\{n,\,1+k(h_A+h_B-2)\}.
\]
\end{proposition}

\begin{proof}
Write $A=x_0+U$ and $B=y_0+V$, and let $W$ be the linear span of all differences of output values.
For $u\in U$ write $\Delta_u^x F(x,y)=F(x+u,y)-F(x,y)$, and define $\Delta_v^y$ similarly for $v\in V$.
Taking $k$ differences in each variable gives
\[
 \Delta_{u_1}^x\cdots\Delta_{u_k}^x
 \Delta_{v_1}^y\cdots\Delta_{v_k}^y P_k(xy)
   =(k!)^2u_1\cdots u_k v_1\cdots v_k .
\]
Indeed the lower powers vanish, and each $k$-th difference of a $k$-th power contributes $k!$ times the product of the increments.
Evaluated at $(x_0,y_0)$, the left-hand side is a sum of output values with coefficients summing to zero. Thus it belongs to $W$.
As $p>k$, $W$ contains the span of all products of $k$ elements of $U$ and $k$ elements of $V$.
Applying Lemma~\ref{lem:products} successively to these $2k$ spaces gives the claimed lower bound for $\dim W$.

For equality take $A=I_{h_A}$ and $B=I_{h_B}$ in the same power basis. Each $P_k(xy)$ is represented by a polynomial in $\theta$ of degree at most $k(h_A+h_B-2)$, which gives the reverse inequality.
\end{proof}

In particular, a nonzero linear coordinate can be constant on two affine sources of dimension $h$ precisely when $2k(h-1)\le n-2$, allowing the coordinate and the sources to be chosen together.
For the example just given one can take the coefficient of $\theta^{n-1}$ as the output. It is identically zero on $I_h\times I_h$.
Thus the limiting rate for a constant coordinate on two equal-dimensional affine sources is $1/(2k)$ when we can choose the coordinate.

\section{The cubic pairing and the paraboloid}\label{app:related}
\subsection{A cubic pairing on the parabola}\label{app:cubic}

For scalar inputs in a prime field, a different pairing of the two parabolas has sharp threshold $1/3$. If we reverse the coordinates of the second parabola, its dot product with the first becomes
\[
 (x,x^2)\cdot(y^2,y)=xy(x+y).
\]
This is also the mixed part of $(x+y)^3$, up to a factor of three.

\begin{corollary}\label{cor:cubicprime}
Let $p>3$ be prime and let $x,y\in\F_p$. Take either the half-interval bit or the least significant bit of $xy(x+y)$ or $(x+y)^3$. Each of the four resulting maps has exponentially small error at every fixed min-entropy rate $\rho>1/3$. At rate $1/3$ each has error bounded below by an absolute constant for all sufficiently large primes $p$.
\end{corollary}

\begin{proof}
Take both inputs from the even integers up to $p^{1/3}/4$. Both polynomials then take even integer values in $[0,p/2)$. Each bit is constant, and Lemma~\ref{lem:enlarge} gives constant error for some pair of sources uniform on sets of size exactly $\lceil p^{1/3}\rceil$.

For the upper bound, reversing the coordinates in Proposition~\ref{prop:bilinear} preserves additive moments and character orthogonality. Thus \eqref{eq:quadratic} gives a power saving for every nonzero frequency of $xy(x+y)$ at rate $\rho>1/3$. The same argument applies to $(x+y)^3=x^3+y^3+3xy(x+y)$. Indeed, in the character sum the factors from $x^3$ and $y^3$ are weights of modulus one. The first application of H\"older removes the $x$-weight. After expanding the even power, the remaining $y$-weights can be removed by taking the absolute value of the inner average over $x$ before the second application of H\"older. The rest of the proof is identical.

These estimates for uniform sources on $N=\lfloor p^\rho\rfloor$ points extend to all sources of rate $\rho$ by convexity. Indeed, the extreme points of the probability vectors with entries at most $1/N$ are uniform on $N$ points. Finally, the residues in the lower half of $\F_p$ form an interval and the even residues are its dilate by two. Their characteristic functions have constant Fourier coefficient $1/2+O(p^{-1})$ and the sum of the absolute values of their normalized Fourier coefficients is $O(\log p)$ by the geometric sum formula. This proves the assertion for both bits.
\end{proof}

\subsection{The paraboloid in three and four ambient dimensions}\label{app:paraboloids}

For comparison, consider the paraboloid in $\F_p^3$. An input is now a point $x\in\F_p^2$, which we lift to $(x,Q(x))$. The dot product of two such points is
\[
 F(x,y)=x\cdot y+Q(x)Q(y),\qquad Q(x)=x_1^2+x_2^2.
\]
If both inputs range over $\{1,\ldots,\lfloor cp^{1/4}\rfloor\}^2$ for small $c>0$, the half-interval bit is constant and each support has size $\asymp\sqrt p$. Allowing the two supports to have different shapes gives a power improvement.

\begin{proposition}\label{prop:paraboloid-prime}
For every sufficiently large odd prime $p$ there are $A,B\subseteq\F_p^2$ with $|A|,|B|\gtrsim p^{2/3}$ such that the half-interval bit and the least significant bit of $F(x,y)$ are each constant on $A\times B$. For either output bit $E$,
\[
 \varepsilon_E(N)\gtrsim pN^{-3/2}\qquad(p^{2/3}\le N\le p^2/2).
\]
For $p^{2/3}\le N\le p$ we have the stronger bound $\varepsilon_E(N)\gtrsim p^{2/3}/N$. In particular there is constant error at rate $1/3$.
\end{proposition}

\begin{proof}
Let $1\le H\le p^{1/3}$ be an integer and take $Y=\{2,4,\dots,2H\}^2$.
Partition $\F_p^2$ according to the intervals containing $[x_1]_p,[x_2]_p$ and $[Q(x)]_p$, of respective lengths $\lfloor cp/H\rfloor,\lfloor cp/H\rfloor$ and $\lfloor cp/H^2\rfloor$, for a sufficiently small absolute $c>0$.
There are $O(H^4)$ sets in the partition. Within any one of them, for fixed $y\in Y$, the even integers
\[
 [x_1]_p y_1+[x_2]_p y_2+[Q(x)]_p(y_1^2+y_2^2)
\]
differ by at most $12cp$. Thus if one residue $[F(x,y)]_p$ is farther than $12cp$ from $0,p/2,p$, both bits are constant throughout the set for this $y$.
For $y\in Y$ the integer $y_1^2+y_2^2$ lies between $8$ and $8H^2<p$ for sufficiently large $p$, so $Q(y)\ne0$. With $x_2$ fixed, the equation $F(x,y)=z$ is quadratic in $x_1$ with nonzero leading coefficient. Each $z$ therefore has at most $2p$ preimages. Taking $c$ small gives at least $3p^2/4$ values of $x$ satisfying the distance condition.
The counting argument in Lemma~\ref{lem:multiplier} now gives a set $A$ of size $\gtrsim p^2/H^4$ for which both bits are constant across $A$ for at least $|Y|/2$ values of $y$. Keeping the most frequent pair of bits gives $B$ of size at least $|Y|/8$.

At $H=\lfloor p^{1/3}\rfloor$ both sets have size $\gtrsim p^{2/3}$. For the error bound take $H=\lfloor(p^2/N)^{1/4}\rfloor$. Then $|A|\gtrsim N$ and $|B|\gtrsim p/\sqrt N$. Apply Lemma~\ref{lem:enlarge}. At the exact rate $1/3$ the prescribed size is $\lceil p^{2/3}\rceil$, since each input has length $2\log_2 p$.

To see the improvement for $N\le p$, suppose first that $p\equiv3\pmod4$. Identify $(x_1,x_2)$ with $x_1+ix_2\in\F_{p^2}$, where $i^2=-1$. Then $Q$ is the field norm. The $p+1$ elements of norm one act freely by multiplication on the nonzero points and preserve both $Q$ and the dot product, by polarization. Rotating both supports by the same element therefore preserves $F$. Discarding zero if necessary, apply the observation following Lemma~\ref{lem:copies} to the rectangle of size $\gtrsim p^{2/3}$ on each side. Every orbit it meets has $p+1$ points. If $p\equiv1\pmod4$, a line of $p$ points on which $F$ vanishes on every pair is described below. This also gives the asserted bound. At $N=\lceil p^{3/4}\rceil$, or rate $3/8$, the resulting error is $\gtrsim p^{-1/12}$.
\end{proof}

If we require the two independent inputs to have the same distribution, we still obtain error at least $1/6$ on a set of size $\gtrsim p^{2/3}$ for either bit. To see this, trim $A$ and $B$ to the same size $M$. If $|A\cap B|\ge M/2$, use the intersection. Otherwise put $A_0=A\setminus B$ and $B_0=B\setminus A$, which are disjoint sets of the same size greater than $M/2$. Write $\beta(S)$ for the signed bias $\Pr(E=1)-1/2$ of two independent uniform samples from $S$. Since $E$ is symmetric and constant on $A_0\times B_0$, we have
\[
 \big|4\beta(A_0\cup B_0)-\beta(A_0)-\beta(B_0)\big|=1.
\]
Thus at least one of the three sets has bias of absolute value at least $1/6$.

There is also an explicit example in which both inputs range over the same set and both bits are constant on every pair. Suppose $p\equiv3\pmod8$, choose $\omega\in\F_p$ with $\omega^2=-2$, and put $R=\lfloor\sqrt p/16\rfloor$. Take
\[
 S=\{(u,\omega v)\in\F_p^2:u,v\in\Z,\ |v|\le R,\ 4R\le u^2-2v^2\le8R\}.
\]
The integer pairs give distinct points of $\F_p^2$, since $|u|\le\sqrt{2R^2+8R}$ and $|v|\le R$. Writing $w=u^2-2v^2$ and $w'=(u')^2-2(v')^2$, the phase on $S\times S$ is represented by
\[
 T=uu'-2vv'+ww'.
\]
The bilinear term has absolute value at most $4R^2+8R$, and $16R^2\le ww'\le64R^2$. Thus $0<T<p/2$. Also $T\equiv uu'+(uu')^2\equiv0\pmod2$, so the least significant bit vanishes on the same set.

The logarithm in its size comes from counting points between the two hyperbolas. For $\sqrt R\le v\le R$, the interval of positive values of $u$ has length
\[
 \sqrt{2v^2+8R}-\sqrt{2v^2+4R}
   =\frac{\sqrt{2}\,R}{v}+O(R^2/v^3).
\]
Rounding the endpoints costs $O(R)$ points in total, as do the values $|v|<\sqrt R$. Summing over $v$ and including both signs of $u,v$ gives
\[
 |S|=2\sqrt2\,R\log R+O(R)
     =\frac{\sqrt p\log p}{8\sqrt2}+O(\sqrt p).
\]
Thus this single explicit support gains $\log_2\log p+O(1)$ bits of entropy over the square.

For $p\equiv3\pmod4$, the half-interval bit has exponentially small error at every rate above $3/8$ by \cite[Remark~7]{Lewko}. Proposition~\ref{prop:paraboloid-prime} leaves open the rates $1/3<\rho\le3/8$ for this bit. When $p\equiv1\pmod4$, the line spanned by $(1,i)$ with $i^2=-1$ gives field output zero on every pair, giving constant error even at rate $1/2$.

Finally, if the inputs themselves lie in $\F_p^3$, the paraboloid lies in $\F_p^4$. Choose $c_1,c_2\in\F_p$ with $c_1^2+c_2^2=-1$. Such a pair exists since the sets $\{z^2:z\in\F_p\}$ and $\{-1-z^2:z\in\F_p\}$ each have $(p+1)/2$ elements. Put $u=(c_1,c_2,1)$ and $v=(c_2,-c_1,0)$, so $u\cdot u=u\cdot v=0$ and $v\cdot v=-1$. Taking $Q(x)=x\cdot x$ as before gives
\[
 F(ru+tv,\,r'u+t'v)=(tt')^2-tt'.
\]
For $1\le t,t'\le p^{1/4}/2$ this is an even integer in $[0,p/2)$. Allowing $r,r'$ to range freely gives one set of size $\gtrsim p^{5/4}$ on whose square both bits are constant. Lemma~\ref{lem:enlarge} gives constant error at the exact rate $5/12$ for inputs of length $3\log_2 p$.


\begin{thebibliography}{99}

\bibitem{AGMR}
O.~Alrabiah, J.~Goodman, J.~Mosheiff and J.~Ribeiro, \emph{Low-degree polynomials are good extractors},
Approximation, Randomization, and Combinatorial Optimization (APPROX/RANDOM 2025), LIPIcs, vol.~353, 2025, pp.~38:1--38:25.

\bibitem{AM}
M.~F.~Atiyah and I.~G.~Macdonald, \emph{Introduction to commutative algebra},
Addison--Wesley, Reading, Mass., 1969.

\bibitem{BSZ}
C.~Bachoc, O.~Serra and G.~Z\'emor, \emph{An analogue of Vosper's theorem for extension fields},
Math. Proc. Cambridge Philos. Soc. \textbf{163} (2017), no.~3, 423--452.

\bibitem{BR}
E.~Ben-Sasson and N.~Ron-Zewi, \emph{From affine to two-source extractors via approximate duality},
SIAM J. Comput. \textbf{44} (2015), no.~6, 1670--1697.

\bibitem{Bombieri}
E.~Bombieri, \emph{On exponential sums in finite fields. II},
Invent. Math. \textbf{47} (1978), no.~1, 29--39.

\bibitem{Bourgain}
J.~Bourgain, \emph{More on the sum-product phenomenon in prime fields and its applications},
Int. J. Number Theory \textbf{1} (2005), no.~1, 1--32.

\bibitem{BDG}
J.~Bourgain, C.~Demeter and L.~Guth, \emph{Proof of the main conjecture in Vinogradov's mean value theorem for degrees higher than three},
Ann. of Math. (2) \textbf{184} (2016), no.~2, 633--682.

\bibitem{BG}
J.~Bourgain and A.~Gamburd, \emph{Uniform expansion bounds for Cayley graphs of $\mathrm{SL}_2(\F_p)$},
Ann. of Math. (2) \textbf{167} (2008), no.~2, 625--642.

\bibitem{Breuillard}
E.~Breuillard, \emph{Approximate subgroups and super-strong approximation},
Groups St Andrews 2013, London Math. Soc. Lecture Note Ser., vol.~422, Cambridge Univ. Press, Cambridge, 2015, pp.~1--50.

\bibitem{ChG}
E.~Chattopadhyay and M.~Gurumukhani, \emph{Low-error two-source extractors via Vinogradov estimates}, to appear in the 18th Innovations in Theoretical Computer Science Conference (ITCS 2027).

\bibitem{CZ}
E.~Chattopadhyay and D.~Zuckerman, \emph{Explicit two-source extractors and resilient functions},
Ann. of Math. (2) \textbf{189} (2019), no.~3, 653--705.

\bibitem{CG}
B.~Chor and O.~Goldreich, \emph{Unbiased bits from sources of weak randomness and probabilistic communication complexity},
SIAM J. Comput. \textbf{17} (1988), no.~2, 230--261.

\bibitem{Dwork}
B.~Dwork, \emph{On the rationality of the zeta function of an algebraic variety},
Amer. J. Math. \textbf{82} (1960), 631--648.

\bibitem{EL}
S.~Eliahou and C.~Lecouvey, \emph{On linear versions of some addition theorems},
Linear Multilinear Algebra \textbf{57} (2009), no.~8, 759--775.

\bibitem{GGMT1}
W.~T.~Gowers, B.~Green, F.~Manners and T.~Tao, \emph{On a conjecture of Marton},
Ann. of Math. (2) \textbf{201} (2025), no.~2, 515--549.

\bibitem{GGMT2}
W.~T.~Gowers, B.~Green, F.~Manners and T.~Tao, \emph{Marton's conjecture in abelian groups with bounded torsion},
Ann. Fac. Sci. Toulouse Math. (6) \textbf{35} (2026), no.~1, 1--33.

\bibitem{HLX}
X.-D.~Hou, K.~H.~Leung and Q.~Xiang, \emph{A generalization of an addition theorem of Kneser},
J. Number Theory \textbf{97} (2002), no.~1, 1--9.

\bibitem{Katz}
N.~M.~Katz, \emph{Sums of Betti numbers in arbitrary characteristic},
Finite Fields Appl. \textbf{7} (2001), no.~1, 29--44.

\bibitem{Lewko49}
M.~Lewko, \emph{An explicit two-source extractor with min-entropy rate near $4/9$},
Mathematika \textbf{65} (2019), no.~4, 950--957.

\bibitem{LewkoRect}
M.~Lewko, \emph{Counting rectangles and an improved restriction estimate for the paraboloid in $\F_p^3$},
Proc. Amer. Math. Soc. \textbf{148} (2020), no.~4, 1535--1543.

\bibitem{Lewko}
M.~Lewko, \emph{A Szemer\'edi--Trotter theorem in arbitrary fields},
preprint, 2026, arXiv:2609.27023.

\bibitem{Li}
X.~Li, \emph{Asymptotically optimal two source extractors for arbitrary constant entropy rate},
preprint, 2026, arXiv:2611.32086.

\bibitem{MM}
S.~Mansfield and A.~Mudgal, \emph{A quadratic Vinogradov mean value theorem in finite fields},
Q. J. Math. \textbf{75} (2024), no.~3, 1007--1029.

\bibitem{PR}
P.~Pudl\'ak and V.~R\"odl, \emph{Pseudorandom sets and explicit constructions of Ramsey graphs},
Complexity of computations and proofs, Quad. Mat., vol.~13, Dept. Math., Seconda Univ. Napoli, Caserta, 2004, pp.~327--346.

\bibitem{RS}
M.~Rudnev and I.~D.~Shkredov, \emph{On the restriction problem for discrete paraboloid in lower dimension},
Adv. Math. \textbf{339} (2018), 657--671.

\bibitem{Varju}
P.~P.~Varj\'u, \emph{Expansion in $\mathrm{SL}_d(\mathcal O_K/I)$, $I$ square-free},
J. Eur. Math. Soc. (JEMS) \textbf{14} (2012), no.~1, 273--305.

\bibitem{WZ}
D.~Wan and D.~Zhang, \emph{Betti number bounds for varieties and exponential sums},
Adv. Math. \textbf{490} (2026), Paper No.~110852.

\bibitem{WooleyCubic}
T.~D.~Wooley, \emph{The cubic case of the main conjecture in Vinogradov's mean value theorem},
Adv. Math. \textbf{294} (2016), 532--561.

\bibitem{WooleyNEC}
T.~D.~Wooley, \emph{Nested efficient congruencing and relatives of Vinogradov's mean value theorem},
Proc. Lond. Math. Soc. (3) \textbf{118} (2019), no.~4, 942--1016.

\end{thebibliography}
\end{document}